\documentclass[10pt]{article}
\usepackage{anyfontsize}
\usepackage{graphicx} 
\usepackage{bbding}
\usepackage{amsfonts}
\usepackage{mathrsfs}
\usepackage{amsthm}
\usepackage{amsmath}

\usepackage{enumerate}
\usepackage{graphicx}
\usepackage{float}
\usepackage{geometry}
\usepackage{amssymb}
\usepackage{bm} 
\usepackage{hyperref}
\usepackage{tikz}
\usepackage{pgfplots}
\pgfplotsset{compat=1.18}
\usepackage{cleveref}
\numberwithin{equation}{section}
\newtheorem{theorem}{Theorem}[section]

\newtheorem{lemma}{Lemma}[section]

\newtheorem{example}{Example}[section]
\newtheorem{proposition}[lemma]{Proposition}

\title{GLOBAL FLOWS FOR KDV HIERARCHY WITH INITIAL CONDITIONS HAVING NO RESTRICTION ON SCHR\"{O}DINGER SPECTRA}
\author{XIANGSHUO KONG, SHINICHI KOTANI, AND SHUO ZHANG }
\date{}

\newcommand{\ip}{\mathrm i}
\newcommand{\dd}{\,\mathrm d}
\newcommand{\la}{\langle}
\newcommand{\ra}{\rangle}
\newcommand{\abs}[1]{\left|#1\right|}
\newcommand{\norm}[1]{\left\lVert#1\right\rVert}
\newcommand{\Oh}{\mathcal O}
\newcommand{\Sh}{\mathcal S}
\newcommand{\Th}{\mathcal T}

\begin{document}
\renewcommand\qedsymbol{$\blacksquare$}
\maketitle

\textbf{Abstract}

When one of the authors constructed the KdV flow for a broad class of initial
data, he imposed a restriction that the underlying Schr\"{o}dinger operators
with the initial data as potentials were bounded from below. We remove this
restriction and establish the flow on a space of initial data that are bounded
from above and possibly unbounded from below.

\textbf{Keywords}: KdV hierarchy; Schr\"{o}dinger operator; Weyl-Titchmarsh
functions; spectral measures; Toeplitz operators; $\tau$-functions; Riccati equation.

\section{Introduction}

In this article the Korteweg-de Vries equation is normalized as%

\begin{equation}
\partial_{t}q=\dfrac{1}{4}\partial_{x}^{3}q-\dfrac{3}{2}q\partial_{x}%
q\text{,}\ \ \ q\left(  x,0\right)  =q_{0}\left(  x\right)  \text{,}%
\ \ \left(  x,t\in\mathbb{R}\right)  \text{.} \label{1.1}%
\end{equation}
Clearly, for any real $c\neq0$ a rescaled function%
\[
u\left(  x,t\right)  =-6c^{2}q\left(  cx,-4c^{3}t\right)
\]
satisfies%
\[
\partial_{t}u+\partial_{x}^{3}u+u\partial_{x}u=0\text{.}%
\]
In \cite{Kotani2024,Kotani2023} it was proved that if the underlying Schr\"{o}dinger operator%
\[
L_{q_{0}}=-\partial_{x}^{2}+q_{0}%
\]
is bounded from below, the KdV equation is uniquely solvable under mild
conditions on the initial $q_{0}$, which admits ergodic initial data
containing almost periodic ones. Recently Zhang \cite{Zhang2026} found that $q_{0}%
\in$ $C_{b}^{n}\left(  \mathbb{R}\right)  $ for a sufficiently large $n\geq1$
can be initial datum. The aim of this article is to remove the condition that
$L_{q}$ is bounded from below, which makes it possible to have a certain class
of functions unbounded from below as initial data.

Menikoff \cite{Menikoff1972} gave a basic idea to study unbounded solutions, and Bondareva \cite{Bondareva1985}
and Bondareva--Shubin \cite{BondarevaShubin1990} obtained an existence theorem on the full time
region $\left(  -\infty,\infty\right)  $ in classes of growing
functions$\ \mathcal{S}^{\alpha}$ with $\alpha\leq1$:
\[
\widetilde{\mathcal{S}}^{\alpha}=\left\{  q_{0}\in C^{\infty}\left(
\mathbb{R}\right)  \text{; \ }q_{0}^{\left(  k\right)  }\left(  x\right)
=o\left(  \left\vert x\right\vert ^{\alpha-k}\right)  \text{ for any }%
k\geq0\right\}  \text{.}%
\]
On the other hand, Kenig, Ponce and Vega \cite{KPV1997} proved further global
well-posedness results on the unilateral time region $[0,\infty)$ for initial
data:%
\[
q_{0}\left(  x\right)  =a_{k}x^{k}+a_{k-1}x^{k-1}+\cdots+a_{1}x+a_{0}+f\left(
x\right)
\]
with $f\in\mathcal{S}\left(  \mathbb{R}\right)  $ (Schwartz class), an odd
integer $k\geq1$ and $a_{k}>0$ (they considered, instead, the equation
$\partial_{t}u+\partial_{x}^{3}u+u\partial_{x}u=0$ for $u\left(  x,t\right)
=-6q\left(  -x,4t\right)  $). However, this solution may explode on the
negative time region as a concrete example%
\[
q_{0}\left(  x\right)  =x\Longrightarrow q\left(  x,t\right)  =\dfrac
{2x}{3t+2}%
\]
implies.

We are interested in global solutions on the full time region $\left(
-\infty,\infty\right)  $. The above previous results suggest the critical
growth order is $o\left(  \left\vert x\right\vert \right)  $ at $x=\pm\infty$.
Furthermore, we would like to remove the decaying condition of the derivatives
of $q_{0}$ in $\widetilde{\mathcal{S}}^{\alpha}$. Our starting point is
\cite{Kotani2024,Kotani2023}, which relates the Sato-Segal-Wilson formalism to the Weyl-Titchmarsh
functions (WT-functions for short) of Schr\"{o}dinger operators. This approach
treats decaying and oscillating potentials in a common setting and includes
substantial classes of bounded smooth ergodic and almost periodic data. Its
use of WT-functions makes it possible to formulate the analytic assumptions
directly in spectral variables.

Since we are not assuming that the associated Schr\"{o}dinger operators
$L_{q_{0}}$ are bounded from below, under the substitution $z=-\zeta^{2}$ of
the spectral parameter $z$, the domain $D_{+}$ should includes $\mathbb{R}\cup
i\mathbb{R}$. We fix a positive odd integer $M$ throught the article. For a
constant $c_{0}>0$ let $\omega\left(  x\right)  $ be a smooth positive even
function on $\mathbb{R}\backslash\left\{  0\right\}  $ satisfying%

\[
\left\vert y\right\vert =\omega\left(  x\right)  =\left\{
\begin{array}
[c]{ll}%
c_{0}+o\left(  1\right)  & \text{as }x\rightarrow\infty\\
O\left(  \left\vert x\right\vert ^{-\left(  M-1\right)  }\right)  & \text{as
}x\rightarrow0
\end{array}
\right.  \text{,}%
\]
and define the domains in $\mathbb{C}$ by%
\[
D_{+}=\left\{  z\in\mathbb{C};\text{ \ }\left\vert \mathfrak{I}z\right\vert
<\omega\left(  \mathfrak{R}z\right)  \right\}  \text{, \ }D_{-}=\left\{
z\in\mathbb{C};\text{ }\left\vert \mathfrak{I}z\right\vert >\omega\left(
\mathfrak{R}z\right)  \right\}  \text{.}%
\]
The boundary of $D_{\pm}$ is denoted by $C$.

\begin{center}
    \begin{tikzpicture}
\begin{axis}[
    axis lines=middle,
    xlabel={$x$},
    ylabel={$y$},
    ticks=none,
    xmin=-5, xmax=5,
    ymin=-5, ymax=5,
    samples=400,
    clip=true,
]

\addplot[blue, thick, domain=-5:-0.05] {abs(x)^(-1/2)};
\addplot[blue, thick, domain=0.05:5]   {abs(x)^(-1/2)};

\addplot[blue, thick, domain=-5:-0.05] {-abs(x)^(-1/2)};
\addplot[blue, thick, domain=0.05:5]   {-abs(x)^(-1/2)};

\node[font=\large] at (axis cs:1,1) {$C$};
\node[font=\large] at (axis cs:0,0) {$D_+$};
\node[font=\large] at (axis cs:2,3) {$D_-$};

\end{axis}
\end{tikzpicture}
\end{center}

\[
\text{the figure of the curve }C
\]
On the Riemann sphere, it forms a closed curve, and in the complex plane, its
interior includes the $x$ and $y$ axes. The curve $C$ is defined so that
$\cosh z^{m}$ remains bounded when $\mathfrak{I}z\rightarrow\infty$ on
$\overline{D}_{+}$ for any positive odd integer $m\leq M$.

Although we follow the framework employed in \cite{Kotani2024,Kotani2023}, multipliers $e^{h(\zeta
)}$ are unbounded on $D_{+}$ and the Hardy spaces and Toeplitz operators must
be modified accordingly. First we define the symbol class for each integer
$N\geq1$%

\[
\mathcal{A}_{N}(C)=\left\{
\begin{array}
[c]{l}%
\boldsymbol{a}\left(  \lambda\right)  =\left(  a_{1}\left(  \lambda\right)
,a_{2}\left(  \lambda\right)  \right)  :a_{j}\text{ are functions on }C\text{
such that for any}\\
c>0\text{ there exists an analytic function }\widetilde{\boldsymbol{a}}\left(
z\right)  \text{ on }D_{+}\text{ satisfying }\\
\sup\limits_{z\in D_{+}}\left\Vert \widetilde{\boldsymbol{a}}\left(  z\right)
\right\Vert +\sup\limits_{\lambda\in C}\left\vert \lambda^{N}\rho_{c}\left(
\lambda\right)  \right\vert \left\Vert \boldsymbol{a}\left(  \lambda\right)
-\widetilde{\boldsymbol{a}}\left(  \lambda\right)  \right\Vert <\infty\text{ }%
\end{array}
\right\}  \text{,}%
\]
where%
\[
\rho_{c}\left(  z\right)  =\left(  K+\cosh z^{M}\right)  ^{c}\text{.}%
\]
The positive constant $K$ is chosen so that $\inf_{z\in D_{+}}$ $\left\vert
\rho_{c}\left(  z\right)  \right\vert >0$. The Hardy spaces $H\left(  D_{\pm
}\right)  $ are defined by%
\[
H\left(  D_{\pm}\right)  =L^{2}\left(  C\right)  \text{-closure of }\left\{
\text{rational functions with no poles on }D_{\pm}\right\}  \text{,}%
\]
and the weighted Hardy space%
\[
H_{N,c}(D_{+})=\left(  z-b\right)  ^{N}\rho_{c}\left(  z\right)
H(D_{+})\text{ with }b\in D_{-}\text{. }%
\]
$z^k\in H_{N,c} (D_+)$ holds for any $0 \leq k \leq N-1 $.
For a rational function $r=q/p$ with polynomials $p$, $q$ set%
\[
\deg r=\deg q-\deg p\text{,}%
\]
and for a polynomial $h$ let $c_{k}\left(  h\right)  $ be the coefficient of
$z^{k}$ of $h$. We introduce semi-groups%
\[%
\begin{tabular}
[c]{l}%
$\Gamma_{m}=\left\{
\begin{array}
[c]{l}%
g=re^{h}\text{; }h\text{ is a real odd polynomial with}\deg h\leq m\text{.
\ }r\text{ is}\\
\text{a rational function such that }r\left(  z\right)  \neq0\text{, }%
\infty\text{ on }\overline{D}_{+}\text{ and}\deg r\leq 0\text{.}%
\end{array}
\right\}  $\\
$\Gamma_{m}^{\operatorname{real}}=\left\{  g\in\Gamma_{m}\text{;
\ }g=\overline{g}\text{ }\right\}  \text{ \ (}\overline{f}\left(  z\right)
=\overline{f\left(  \overline{z}\right)  }\text{)}$%
\end{tabular}
\]
For $\boldsymbol{a}\in\mathcal{A}_{N}(C)$ and $g=re^{h}\in\Gamma_{m}$ with
$m\leq M$ define a Toeplitz operator $T(g\boldsymbol{a})$ on $H_{N,c}(D_{+})$ by%

\[
T(g\boldsymbol{a})u=\mathfrak{p}_{+}(g(\boldsymbol{a}-\widetilde
{\boldsymbol{a}})u)+g\widetilde{\boldsymbol{a}}u\text{,}%
\]
where $\mathfrak{p}_{+}$ denotes the projection from $L^{2}\left(  C\right)  $
to $H\left(  D_{+}\right)  $ and $\boldsymbol{a}u=a_{1}u_{e}+a_{2}u_{o}$ with
$u_{e}\left(  \lambda\right)  =\left(  u\left(  \lambda\right)  +u\left(
-\lambda\right)  \right)  /2$, $u_{o}=u-u_{e}$. The first term of the right
hand belongs to $H(D_{+})$. The operator $T(g\boldsymbol{a})$ does not depend
on the choice of the compensator $\widetilde{\boldsymbol{a}}$ and satisfies%
\[
g^{-1}T(g\boldsymbol{a}):H_{N,c}(D_{+})\rightarrow H_{N,c}(D_{+})\text{ \ if
}-c_{m}\left(  h\right)  \leq c\text{,}%
\]
since, for $g=re^{h}\in\Gamma_{m}$ with $m\leq M$ it holds that $gu\in
H_{N,c}(D_{+})$ for $u\in H(D_{+})$ if $\deg r\leq N$, $c_{m}\left(  h\right)
\leq c$. Define%
\[
\mathcal{A}_{N}^{inv}(C)=\left\{  \boldsymbol{a}\in\mathcal{A}_{N}(C)\text{;
}T(\boldsymbol{a})\text{ is invertible on }H_{N,c}(D_{+})\text{ for any
}c>0\right\}  \text{.}%
\]
An important subclass $\mathcal{B}_{N}$ of $\mathcal{A}_{N}^{inv}(C)$ is given
in Section 2.

For any rational function $r$ of $\Gamma_{m}$ the operator $r^{-1}%
T(r\boldsymbol{a})T(\boldsymbol{a})^{-1}-I$ is of finite rank in
$H_{N,c}(D_{+})$ if $\boldsymbol{a}\in\mathcal{A}_{N}^{inv}(C)$, hence we can
define the tau-function:%
\[
\tau_{\boldsymbol{a}}\left(  r\right)  =\det\left(  r^{-1}T(r\boldsymbol{a}%
)T(\boldsymbol{a})^{-1}\right)  \text{ for }\boldsymbol{a}\in\mathcal{A}%
_{N}^{inv}(C)\text{.}%
\]
In order to show that $g^{-1}T(g\boldsymbol{a})$ is invertible we introduce%
\[
\mathcal{A}_{N,+}^{inv}(C)=\left\{
\begin{array}
[c]{l}%
\boldsymbol{a}\in\mathcal{A}_{N}^{inv}(C)\text{; }\boldsymbol{a}%
=\overline{\boldsymbol{a}}\text{ and }\tau_{\boldsymbol{a}}\left(  r\right)
\geq0\text{ for any real rational }\\
\text{function }r\text{ with no poles nor zeros in }\overline{D}_{+}%
\end{array}
\right\}  \text{.}%
\]
In Section 4 we prove that $g^{-1}T(g\boldsymbol{a})$ is invertible for
$\boldsymbol{a}\in\mathcal{A}_{N,+}^{inv}(C)$ and $g\in\Gamma_{M}%
^{\operatorname{real}}$. In this case we define
\[
\varphi_{g\boldsymbol{a}}(z)=g(z)\boldsymbol{a}(z)(g^{-1}T(g\boldsymbol{a}%
))^{-1}(g^{-1}\cdot1)(z)-1\in H(D_{-})\text{,}%
\]
and%
\[
\kappa_{1}\left(  g\boldsymbol{a}\right)  =\lim_{D_{-}\ni z\rightarrow\infty
}z\varphi_{g\boldsymbol{a}}(z)\text{.}%
\]
The initial space for the KdV flow is given by%
\[
\mathcal{Q}_{N}^{\left(  M\right)  }=\left\{  q\text{; \ }q\left(  x\right)
=-2\partial_{x}\kappa_{1}\left(  e_{x}\boldsymbol{a}\right)  \text{ for an
}\boldsymbol{a}\in\mathcal{A}_{N,+}^{inv}(C)\right\}  \text{,}%
\]
where $e_{x}\left(  z\right)  =e^{xz}$. The odd integer $M$ and integer
$N\geq 3$ describe the degree of the KdV hierachy and the smoothness of $q$
respectively. Our first theorem is

\begin{theorem}
\label{Thm: theorem1}The followings are valid.\newline(i) $\mathcal{Q}_{N}^{\left(
M\right)  }\subset C^{N-3}\left(  \mathbb{R}\right)  $ holds. \newline(ii) For
$q\in\mathcal{Q}_{N}^{\left(  M\right)  }$ and $g\in\Gamma_{M}%
^{\operatorname{real}}\boldsymbol{\ }$define%
\[
\left(  K\left(  g\right)  q\right)  \left(  x\right)  =-2\partial_{x}%
\kappa_{1}\left(  e_{x}g\boldsymbol{a}\right)  \text{ for }q\left(  x\right)
=-2\partial_{x}\kappa_{1}\left(  e_{x}\boldsymbol{a}\right)  \text{.}%
\]
For any odd integer $m\leq M$ with $m+3<N$ the function $\left(  K\left(  e^{tz^{m}%
}\right)  q\right)  \left(  x\right)  $ is $C^{1}$ in $t$ and $C^{N-3}$ in
$x$, and it satisfies the $\left(  m+1\right)  /2$ th KdV equation for
$q\in\mathcal{Q}_{N}^{\left(  M\right)  }$. Especially, $\left(  K\left(
e^{tz^{3}}\right)  q\right)  \left(  x\right)  $ is a solution to the KdV
equation (\ref{1.1}) for $q\in\mathcal{Q}_{N}^{\left(  3\right)  }$ if
$N\geq6$.
\end{theorem}

To give concrete examples of elements of $\mathcal{Q}_{N}^{\left(  M\right)
}$ we introduce $m$-functions. For $\boldsymbol{a}\in\mathcal{A}_{N,+}%
^{inv}(C)$, $g\in\Gamma_{M}^{\operatorname{real}}$ let%
\[
\psi_{g\boldsymbol{a}}(z)=g(z)\boldsymbol{a}(z)(g^{-1}T(g\boldsymbol{a}%
))^{-1}(g^{-1}\cdot z)(z)-z\in H(D_{-})\text{,}%
\]
and set%
\[
m_{g\boldsymbol{a}}(z)=\dfrac{z+\psi_{g\boldsymbol{a}}(z)}{1+\varphi
_{g\boldsymbol{a}}(z)}+\kappa_{1}\left(  g\boldsymbol{a}\right)  \text{.}%
\]
We call this $m_{g\boldsymbol{a}}$ as the $m$-function of $g\boldsymbol{a}$.
In Section 6 it will be proved that $m_{\boldsymbol{a}}$ is descrived by the
WT-functions $m_{\pm}$ of the Schr\"{o}dinger operator $L_{q}$ with $q\left(
x\right)  =-2\partial_{x}\kappa_{1}\left(  e_{x}\boldsymbol{a}\right)  $ as
follows:%
\begin{equation}
m_{\boldsymbol{a}}(z)=\left\{
\begin{array}
[c]{cc}%
-m_{+}\left(  -z^{2}\right)  & \text{if \ }\mathfrak{R}z>0\\
m_{-}\left(  -z^{2}\right)  & \text{if \ }\mathfrak{R}z<0
\end{array}
\right.  \text{.} \label{1.2}%
\end{equation}
Moreover, the symbol $\boldsymbol{b}\left(  \lambda\right)  =\left(
1,m_{\boldsymbol{a}}(\lambda)/\lambda\right)  $ gives the same $q$ (see
Section 6). Therefore, a sufficient condition for $q$ to be an element of
$\mathcal{Q}_{N}^{\left(  M\right)  }$ can be stated by its WT-functions. For
an analytic function $m\left(  z\right)  $ on $\mathbb{C}\backslash\left(
\mathbb{R}\cup i\mathbb{R}\right)  $ consider the following conditions:

\noindent(M.1) \ $m\left(  z\right)  =\overline{m\left(  \overline{z}\right)
}$ \ and \ $\mathfrak{I\ }m\left(  z\right)  /\mathfrak{I\ }z>0$ on
$\mathbb{C}\backslash\left(  \mathbb{R}\cup i\mathbb{R}\right)  $%
.\newline\noindent(M.2) \ $m\left(  z\right)  =z+\sum_{1\leq k\leq N-2}%
m_{k}z^{-k}+O\left(  z^{-N+1}\right)  $ on every $D_{-}$.\newline(M.3) \ the
spectral measures $\sigma_{\pm}$ of the Herglotz functions $m_{\pm}\left(
z\right)  =\mp m\left(  \sqrt{-z}\right)  $ satisfies%
\[
\int_{-\infty}^{0}e^{c\left\vert \lambda\right\vert ^{M/2}}\sigma_{\pm}\left(
d\lambda\right)  <\infty\text{ for any }c>0\text{.}%
\]

\begin{theorem}
\label{Thm: theorem2}For the WT-functions $m_{\pm}$ of a potential $q$ define $m\left(
z\right)  $ by the right hand of (\ref{1.2}). If $m\left(  z\right)  $
satisfies the conditions (M.1)\symbol{126}(M.3), then $q\in\mathcal{Q}%
_{N}^{\left(  M\right)  }$. Moreover, $m_{K(g)q}$ satisfies (M.1), (M.3) and(M.2) for $N-\frac{M+1}{2}$.
\end{theorem}

A sufficient condition for the WT-functions satisfying the conditions (M.1),
(M.2), (M.3) in terms of $q$ is given by

\begin{theorem}
\label{Thm: theorem3}For $q\in C^{k}\left(  \mathbb{R}\right)  $ with even integer $k\geq M+N+2$
assume the following conditions:\newline(A) $C_{0}:=\sup_{x\in\mathbb{R}%
}q\left(  x\right)  <\infty$, and there exists a constant $c>0$ and
$\alpha<2/\left(  M-1\right)  $ such that $q\left(  x\right)  \geq-c\left(
1+\left\vert x\right\vert \right)  ^{\alpha}$ holds for every sufficiently
large $\left\vert x\right\vert $.\newline(B) Set $Q\left(  x\right)
=1+C_{0}-q\left(  x\right)  \left(  \geq1\right)  $ and $a=\left(  1-\left(
M+N+2\right)  /k\right)  /2\left(  \in\lbrack0,1/2)\right)  $.%
\[
\sup_{x\in\mathbb{R}}\dfrac{\left\vert q^{\left(  j\right)  }\left(  x\right)
\right\vert }{Q\left(  x\right)  ^{1+ja}}<\infty\text{ for }j=1\text{,
}2\text{,}\cdots\text{, }k\text{. }%
\]
Then, $q\in\mathcal{Q}_{N}^{\left(  M\right)  }$ holds, since the associated
$m$-function satisfies (M.1), (M.2), (M.3).
\end{theorem}

The condition (M.2) is proved under (B) by applying the approach employed by
Zhang \cite{Zhang2026}, whereas (M.3) is obtained from the estimate of
exponential integrals with respect to the spectral measure on $\left(
-\infty,0\right)  $ in terms of the energy functional, which was established
by Gesztesy-Simon \cite{GesztesySimon2000}.

The condition (B) includes many initial data $q$. Typical examples are
provided below.

\begin{example}
A priori assume $q\in C^{k}\left(  \mathbb{R}\right)  $.\newline1. Bounded
smooth potentials $C_{b}^{k}\left(  \mathbb{R}\right)  $ (Zhang \cite{Zhang2026}'s case)\newline2. Potentials decaying to $-\infty$: For
$0<\alpha<2/\left(  M-1\right)  $, $1\leq\beta<1+\alpha/2$%
\[
q\left(  x\right)  =-\left\{
\begin{array}
[c]{ll}%
x^{\alpha}u_{+}\left(  x^{\beta}\right)  & \text{if }x\geq1\\
\left(  -x\right)  ^{\alpha}u_{-}\left(  \left(  -x\right)  ^{\beta}\right)  &
\text{if }x\leq-1
\end{array}
\right.  +v\left(  x\right)
\]
with $u_{\pm}$, $v\in C_{b}^{k}\left(  \mathbb{R}\right)  $ satisfying
$\inf_{x\geq1}u_{\pm}\left(  x\right)  >0$. The derivatives of $q$ may diverge
as $x\rightarrow\pm\infty$. \newline3. Sparse negative wells with unbounded
depths:%
\[
q\left(  x\right)  =-\sum_{n\geq2}nw\left(  \dfrac{x-2^{n^{2}}}{n}\right)
\text{ \ with }0\leq w\in C_{c}^{k}\left(  \left(  -1,1\right)  \right)
\text{.}%
\]

\end{example}

Killip-Murphy-Visan \cite{Killip2020} constructed a solution to the KdV equation with
Gaussian white noise initial condition in 2020. We are interested in this
contruction from the point of view of Sato-Segal-Wilson theory. In this
respect we are expecting that the conditions (M.1)-(M.3) are verified for the
Gaussian white noise potential. (M.1) trivially holds. (M.3) is valid for
$M=3$ and small $c>0$, which follows from the explicit formula for the density
of states. Therefore, the main difficulty lies in (M.2) for $N=1$, or more
precisely, the possibility of defining $\kappa_{1}\left(  e_{x,t}%
\boldsymbol{a}\right)  $ with $e_{x,t}\left(  z\right)  =e^{xz+tz^{3}}$.

The paper is organized as follows. Section 2 introduces the weighted Hardy
spaces and Toeplitz operators. Sections 3 and 4 develop the regularized
tau-functions and their positivity and nondegeneracy properties. Section 5
derives the Schr\"{o}dinger equation and the KdV evolution from the
reconstruction. Up to this section we follow the framework taken by [8,9]
and minimize the proofs. Section 6 relates the construction of the flow to WT-
functions and spectral measures and proves the main results on the flow.
Section 7 establishes the thin-region asymptotics and the negative spectral
moment estimates for concrete classes of potentials.

\section{Hardy spaces and Toeplitz operators}
We introduce the weighted Hardy spaces and their Toeplitz operators, and give conditions on $g_1,g_2$ under which $g_1^{-1}T(g_1a)-g_2^{-1}T(g_2a)$ is Hilbert--Schmidt.

The Hardy spaces associated with $C$ are given by
\begin{equation}
    H_\pm=H(D_{\pm})=\text{ the closure in }L^2(C)\text{ of rational function with no poles in }D_\pm,
\end{equation}
and $L^2(C)=H_+\oplus H_-(\text{direct sum})$ with $\mathfrak{p}_\pm$ Cauchy projections from $L^2(C)$ to $H_\pm$  given by 
\begin{equation}
    \mathfrak{p}_{\pm}f=\pm\frac{1}{2\pi i}\int_{C}\frac{f(\lambda)}{\lambda-z}\dd\lambda,\qquad z\in D_\pm
\end{equation}

To accommodate the growth of $e^h$ on $C$ and the fact that $z^n\notin H(D_+)$, we work in the weighted Hilbert space $\lambda^N\rho_c(\lambda)L^2(C)$. With $\rho_c(z)=(\cosh z^M+K)^c$, define
\begin{equation}
    H_{N,c}(D_+)=(z-b)^N\rho_c(z)H(D_+),\qquad b\in D_-
\end{equation}
where $K=K(M)>0$ is chosen so that $\left(\cosh{z^M+K}\right)^c$ has no zeros in $D_+$ and 
\begin{equation}
    z^k\in H_{N,c}(D_+), \qquad\text{for any integer }k\leq N-1.
\end{equation}

For the symbols on the weighted Hardy space, define the class 
\begin{equation}
   \mathcal{A}_{N}(C)=\left\{
    \begin{aligned}
        &\bm{a}(\lambda):\bm{a}\text{ is a bounded function on }C,\text{ and for any }c>0,\text{ there exists an}\\
        &\text{analytic function }\tilde{\bm{a}}\text{ satisfying }\sup_{\lambda\in D_+}|\tilde{\bm{a}}(z)|<\infty,\ \sup_{\lambda\in C}|\lambda^N\rho_c(\lambda)||\bm{a}(\lambda)-\tilde{\bm{a}}(\lambda)|<\infty
    \end{aligned}
    \right\}
\end{equation}
associated with $H_{N,c,+}$. We equip $H_{N,c}(D_+)$ with the norm

\begin{equation}
    \|f\|_{H_{N,c}(D_+)}=\sqrt{\int_C\big|f(z)(z-b)^{-N}\rho_c(z)^{-1}\big|^2|dz|}.
\end{equation}

\begin{lemma}
    $H(D_+)$ is dense in $H_{N,c}(D_+)$.
\end{lemma}
\begin{proof}
    Fix $u\in H_{N,c}(D_+)$ and set $u_n(z)=\frac{u(z)}{(z-b)^N\rho_c(z)}f_n(z)$, where $f_n(z)=e^{\frac{-z^{4M}}{n}}(z-b)^N\rho_c(z)$. Since $f_n(z)$ is bounded and analytic in $D_+$, we have $u_n\in H(D_+)$. Moreover, $u_n(z)\rightarrow u(z)$ pointwise in $D_+$. Thus 
    \begin{equation}
        \|u_n-u\|_{H_{N,+}(D_+)}=\int_{C}|u(z)-u_n(z)|^2\frac{|dz|}{|z^N\rho_c(z)|^2}=\int_C\abs{\frac{(z-b)^N}{z^N}}\frac{|u(z)|^2}{|\rho_c(z)|^2}\left|e^{\frac{-z^{4M}}{n}}-1\right|^2|dz|\rightarrow0
    \end{equation}
    by the dominated convergence theorem.
\end{proof}

We form the direct sum
\begin{equation}
    H_{N,c}=H_{N,c}(D_+)\bigoplus H(D_-)
\end{equation}
and extend the projections to $H_{N,c}$ by setting
\begin{equation}
    \mathfrak{p}_+f=f_+,\mathfrak{p}_-f=f_- \text{ if } f=f_++f_-\text{ with } f_+\in H_{N,c}(D_+),f_-\in H(D_-).
\end{equation}
The norm on the direct sum, with its positive component in $H_{N,c}(D_+)$, is
\begin{equation}
    \|f\|_{H_{N,c}}=\sqrt{\int_C|f_+(\lambda)|^2|\lambda^N\rho_{c}(\lambda)|^{-2}|\dd\lambda|}+\sqrt{\int_{C}|f_-(\lambda)|^2|\dd \lambda|}
\end{equation}

We shall use the family
\begin{equation}
    \Gamma_m=\left\{
    \begin{aligned}
        &g=re^h;r\text{ is a rational function such that}\\
        &(\romannumeral1)\ r\text{ has no poles nor zeros on }\mathbb{R}\cup i\mathbb{R}\\
        &(\romannumeral2)\ \deg{r}\leq0\\
        &(\romannumeral3)\ h \text{ is real polynomials of degree}\leq m
    \end{aligned}
    \right\}
\end{equation}
Let $g\in \Gamma_m$ and $\bm{a}\in \mathcal{A}_{N}(C)$, with $m\leq M$. The Toeplitz operator $T(g\bm{a})$ on $H_{N,c}$ is defined by
\begin{equation}
   T(g\bm{a})u=\mathfrak{p}_+(g\bm{a}u)=\mathfrak{p}_+(g(\bm{a}-\tilde{\bm{a}})u)+g\tilde{\bm{a}}u\in gH_{N,c}(D_+).
\end{equation}
Here $\bm{a}u=a_1u_e+a_2u_o$.

Splitting the Cauchy kernel gives
\begin{equation}\label{eq:calculation of T(ga)-gT(a)}
\begin{aligned}
    T(g\bm{a})u(z)&=\frac{1}{2\pi i}\int_{C}\frac{g(\lambda)}{\lambda-z}(\bm{a}(\lambda)-\tilde{\bm{a}}(\lambda))u(\lambda)d\lambda+g\tilde{\bm{a}}u(z)\\
    &=\frac{1}{2\pi i}\int_{C}\frac{g(\lambda)-g(z)}{\lambda-z}(\bm{a}(\lambda)-\tilde{\bm{a}}(\lambda))u(\lambda)d\lambda+g\mathfrak{p}_+(\bm{a}-\tilde{\bm{a}}u)(z)+g\tilde{\bm{a}}u(z)\\
    &=\frac{1}{2\pi i}\int_{C}\frac{g(\lambda)-g(z)}{\lambda-z}(\bm{a}(\lambda)-\tilde{\bm{a}}(\lambda))u(\lambda)d\lambda+gT(\bm{a})u(z).
\end{aligned}
\end{equation}

\begin{lemma}
    Let $g_1=r_1e^{h_1}, g_2=r_2e^{h_2}\in\Gamma_m$, where the rational factors have expansions $r_j=\sum_{i=n_j}^{-\infty}a_{j,i}z^{i}$ and the odd polynomials are given by $h_j=\sum_{i=0}^{m_j}b_{j,i} z^i$. Set $n=\min\{n_1,n_2\}\leq0$ and $C=\max_{j}\{|b_{j,m_j}|\}$. Fix $\varepsilon>0$ such that $|e^{h_j(z)}|\lesssim e^{(C+\varepsilon)|\Re z|^m} $. If $M,N,c$ satisfy
    \begin{equation}\label{eq:condition of MNc}
        N\geq\max\{m,-n+1\},\qquad M\geq m,\qquad c>C+\varepsilon\text{ when }M=m
    \end{equation}  then $g_1^{-1}T(g_1\bm{a})-g_2^{-1}T(g_2\bm{a})$ is Hilbert--Schmidt.
\end{lemma}
\begin{proof}

For $g_1=r_1e^{h_1}, g_2=r_2e^{h_2}\in\Gamma_m$ with $r_j=\sum_{i=n_j}^{-\infty}a_{j,i}z^{i}, h_j=\sum_{i=0}^{m_j}b_{j,i} z^i$, the preceding identity gives 
\begin{equation}\label{2.10}
\begin{aligned}
    &g_1^{-1}T(g_1a)u(z)-g_2^{-1}T(g_2a)u(z)\\
    &=\frac{1}{2\pi i}\int_{C}\frac{g_1^{-1}(z)g_1(\lambda)-g_2^{-1}(z)g_2(\lambda)}{\lambda-z}(a(\lambda)-\tilde{a}(\lambda))u(\lambda)\dd\lambda\\
    &=\frac{1}{2\pi i}\int_{C}\frac{g_1^{-1}(z)g_1(\lambda)-g_2^{-1}(z)g_2(\lambda)}{\lambda-z}(a(\lambda)-\tilde{a}(\lambda))\lambda^{2N}\rho^2_c(\lambda)u(\lambda)\lambda^{-2N}\rho^{-2}_c(\lambda)\dd\lambda
\end{aligned}
\end{equation}
 The function $\cosh{z^M+K}$ remains bounded as $\Im{z}\rightarrow\infty$ along $C$. Using the weight estimate $|z^N\rho_c(z)|=|z^{N}\left(\cosh{(z^M+K)}\right)^c|\sim |z|^N e^{c|\Re{z}|^M}$, we obtain
\begin{equation}\label{eq:HS norm}
\begin{aligned}
    &\|g_1^{-1}T(g_1a)u(z)-g_2^{-1}T(g_2a)u(z)\|^2_{H.S.}\\
    &\leq\frac{1}{4\pi^2}\int_{C^2}\left|\frac{g_1^{-1}(z)g_1(\lambda)-g_2^{-1}(z)g_2(\lambda)}{\lambda-z}(a(\lambda)-\tilde{a}(\lambda))\lambda^{2N}\rho_c^2(\lambda)\right|^2|z^N\rho_c(z)|^{-2}|\lambda^N\rho_c(\lambda)|^{-2}|dz||d\lambda|\\
    &\leq c_a\frac{1}{4\pi^2}\int_{C^2}\left|\frac{g_1^{-1}(z)g_1(\lambda)-g_2^{-1}(z)g_2(\lambda)}{\lambda-z}\right|^2|z^{N}\rho_c(z)|^{-2}|\lambda|^{2N}|\lambda^N\rho_c(\lambda)|^{-2}|dz||d\lambda|
\end{aligned}    
\end{equation}
since $\tilde{a}$ can be chosen so that $\rho_{2c}(a-\tilde{a})$ is bounded on $C$. The weight in the last integral satisfies 
\begin{equation}
    |z^N\rho_c(z)|^{-2}|\lambda|^{2N}|\lambda^N\rho_c(\lambda)|^{-2}\sim |z|^{-2N}e^{-2c\left(|\Re{z}|^M+|\Re{\lambda}|^M\right)}
\end{equation}
as $z,\lambda\rightarrow\infty$.

Subtracting the diagonal values gives the identity 
\begin{equation}
    \frac{g_1^{-1}(z)g_1(\lambda)-g_2^{-1}(z)g_2(\lambda)}{\lambda-z}=\frac{g_1^{-1}(z)(g_1(\lambda)-g_1(z))-g_2^{-1}(z)(g_2(\lambda)-g_2(z))}{\lambda-z}.
\end{equation}
Split $C^2$ into the two regions $C^2=\{(z,\lambda):|z-\lambda|\leq\varepsilon|\lambda|\}\cup\{(z,\lambda):|z-\lambda|>\varepsilon|\lambda|\}$. Let $I_1,I_2$ denote the contributions to the integral in \eqref{eq:HS norm} from the respective regions of $C^2$. Define the metric on $\Gamma_m$ by 
\begin{equation}\begin{aligned}
    d_c(g_1,g_2)&=\sqrt{\frac{1}{4\pi^2}\int_{C^2}\left|\frac{g_1^{-1}(z)g_1(\lambda)-g_2^{-1}(z)g_2(\lambda)}{\lambda-z}\right|^2|z^N\rho_c(z)|^{-2}|\lambda|^{2N}|\lambda^N\rho_c(\lambda)|^{-2}|dz||d\lambda|}\\
    &\leq\sqrt{I_1+I_2}.
\end{aligned}
\end{equation}

On the first region, use the identity $\frac{g_j(\lambda)-g_j(z)}{\lambda-z}=\int_{0}^{1}g_j'((\lambda-z)t+z)dt$. With $C=\max_j\{|b_{j,m_j}|\}$, we have \begin{equation}
    g_j'(z)=(r_j'(z)+r_j(z)h_j'(z))e^{h_j(z)}\leq C_{\varepsilon}z^{n_j+m-1}e^{(C+\varepsilon)|\Re z|^m}
\end{equation} and 
\begin{equation}\begin{aligned}
    \Delta_1&\leq\frac{1}{4\pi^2}\int_{|z-\lambda|<\epsilon|\lambda|}\left|\frac{g_1^{-1}(z)g_1(\lambda)-g_2^{-1}(z)g_2(\lambda)}{\lambda-z}\right|^2|z^N\rho_c(z)|^{-2}|\lambda|^{2N}|\lambda^N\rho_c(\lambda)|^{-2}|dz||d\lambda|\\
    &\lesssim\frac{1}{4\pi^2}\int_{|z-\lambda|<\varepsilon|\lambda|}\left|z^{-n_j}e^{-h_j(z)}z^{n_j+m-1}e^{C|\Re{z}|^m}\right|^2\left|z^{-N}e^{-2c|\Re{z}|^M}\right|^2|dz||d\lambda|\\
    &\lesssim\frac{1}{4\pi^2}\int_{|z-\lambda|<\varepsilon|\lambda|}\left|z^{-N+m-1}e^{-2c|\Re{z}|^M+2C|\Re{z}|^m}\right|^2|dz||d\lambda|
\end{aligned}\end{equation}
Thus $\Delta_1<\infty $ provided that $N>m$ and $M\geq m$, with $c>C$ when $M=m$.

On the second region, where $|z-\lambda|>\epsilon|\lambda|$, we have 
\begin{equation}
    \left|\frac{g_1^{-1}(z)g_1(\lambda)-g_2^{-1}(z)g_2(\lambda)}{\lambda-z}\right|\leq c_1|\lambda|^{-1}(|z|^{-n}|\lambda|^{n}e^{C(|\Re z|^m+|\Re\lambda|^m)})
\end{equation}
where $n=\min\{n_1,n_2\}\leq0$.

\begin{equation}\begin{aligned}
    \Delta_2&\leq\frac{1}{4\pi^2}\int_{|z-\lambda|>\epsilon|\lambda|}\left|\frac{g_1^{-1}(z)g_1(\lambda)-g_2^{-1}(z)g_2(\lambda)}{\lambda-z}\right|^2|z^N\rho_c(z)|^{-2}|\lambda|^{2N}|\lambda^N\rho_c(\lambda)|^{-2}|dz||d\lambda|\\
    &\lesssim\frac{1}{4\pi^2}\int_{|z-\lambda|>\epsilon|\lambda|}\left(|\lambda|^{-2}|z|^{-2n}|\lambda|^{2n}e^{2C(|\Re z|^m+|\Re\lambda|^m)}\right)|z|^{-2N} e^{-2c(|\Re{z}|^M+|\Re{\lambda}|^M)}|dz||d\lambda|\\
    &=\frac{1}{4\pi^2}\int_{|z-\lambda|>\epsilon|\lambda|}|\lambda|^{2(-1+n)}|z|^{2(-n-N)}e^{2C(|\Re z|^m+|\Re\lambda|^m)-2c(|\Re{z}|^M+|\Re{\lambda}|^M)}|dz||d\lambda|.
\end{aligned}
\end{equation}
Hence $\Delta_2<\infty$ provided that $N\geq-n+1,M\geq m$, with $c>C$ when $M=m$.
\end{proof}
Within $\mathcal{A}^{inv}_{N}$, consider the class 
\begin{equation}
    \mathcal{B}_{N}=\left\{\begin{aligned}
        \bm{m}(z)=(m_1(z),m_2(z)):\bm{m}\text{ is holomorphic on }\mathbb{C}\setminus(\mathbb{R}\cup i\mathbb{R})\text{ and satisfies (i)(ii)}
    \end{aligned}\right\}
\end{equation}
where conditions (i) and (ii) are as follows.\\
(i) $\bm{m(z)}=\bm{1}+\sum_{1\leq k< N}\bm{m}_kz^{-k}+\Oh(z^{-N})$ in $D_-$, where $\bm{1}=(1,1),\bm{m}_k\in\mathbb{C}^2$.\\
(ii) $m_1(z)m_2(-z)+m_1(-z)m_2(z)\neq0$ on $\mathbb{C}\setminus(\mathbb{R}\cup i\mathbb{R})$.

For $\bm{m},\bm{n}$ in this class, define their product and the element $\hat{\bm m}$ by
\begin{equation}
    \left\{\begin{aligned}
        (\bm{m}\cdot\bm{n})(z)=(m_1(z)n_{1,e}(z)+m_2(z)n_{1,o}(z),m_1(z)n_{2,o}(z)+m_2(z)n_{2,e}(z))\\
        \hat{\bm{m}}(z)=\left(\frac{2(m_{2,e}(z)-m_{1,o}(z))}{m_1(z)m_2(-z)+m_1(-z)m_2(z)},\frac{2(m_{1,e}(z)-m_{2,o}(z))}{m_1(z)m_2(-z)+m_1(-z)m_2(z)}\right)
    \end{aligned}\right.
\end{equation}
\begin{lemma}
    Under this product, $\mathcal{B}_{N}$ is a group, and the associated operators $\{T(\bm{m}):\bm{m}\in\mathcal{B}_{N}\}$ satisfy
    \begin{equation}
        T(\bm{m}\cdot\bm{n})=T(\bm{m})T(\bm{n}).
    \end{equation}
    
\end{lemma}
\begin{proof}
The even--odd product rule gives
$\bm{m}\cdot\hat{\bm{m}}=\hat{\bm{m}}\cdot\bm{m}=(1,1)$, as in
\cite[Lemma~7]{Kotani2023}. For the operator identity, the even and odd
parts of $\mathfrak{p}_-(\bm{n}u)$ belong to $H(D_-)$. Multiplication
by $m_1$ or $m_2$ preserves this space, so
\[
 T(\bm{m}\cdot\bm{n})u
 =\mathfrak{p}_+(\bm{m}(\bm{n}u))
 =\mathfrak{p}_+(\bm{m}\,\mathfrak{p}_+(\bm{n}u))
 =T(\bm{m})T(\bm{n})u .
\]
\end{proof}

\section{Tau-functions}
If $g^{-1}T(g\bm{a})T(\bm{a})^{-1}-I$ is trace class, the $\tau $-function is given by the Fredholm determinant
\begin{equation}\label{eq:def of tau}
    \tau_{\bm{a}}(g)=\det(g^{-1}T(g\bm{a})T(\bm{a})^{-1})
\end{equation}
For the Hilbert--Schmidt difference $g_1^{-1}T(g_1a)-g_2^{-1}T(g_2a)$, we use the regularized determinant $\det_2$, given by
\begin{equation}
    {\det}_2(I+A)=\det((I+A)e^{-A}) .
\end{equation}
In particular, if $A$ is trace class, then
\begin{equation}
    {\det}_2(I+A)=\det(I+A)e^{-\operatorname{tr}A} .
\end{equation}
For a Hilbert--Schmidt operator $A$, $I+A$ is invertible if and only if $\det_2(I+A)\neq 0$. Thus $\det_2$ gives a criterion for the existence of $T(ga)^{-1}$. We define
\begin{equation}\label{eq:def of tau2}
    \tau_{\bm{a}}^{(2)}(g)={\det}_2(g^{-1}T(g\bm{a})T(a)^{-1})
\end{equation}
for $\bm{a}\in\mathcal{A}^{inv}_{N},g=re^h\in\Gamma_m$.

For the operator $T(g\bm{a})$ and its $\tau$-functions, we compare the choices of $M,N,c$ using  of \cite[Lemma~4.8]{Kotani2024}.
\begin{lemma}\label{lem:subspace}
    Let $H,H_1$ be Hilbert spaces and suppose that $H_1$ is a vector subspace of $H$ and that $H_1$ is dense in $H$. Let $A$ be a Hilbert--Schmidt operator on $H$ satisfying 
    \begin{equation}
        AH_1\subset H_1 \text{ and } A_{H_1}=A|_{H_1}\text{ is of Hilbert-Schmidt class in } H_1
    \end{equation}
    Then
    \begin{equation}
        {\det}_2(I_H+A)={\det}_2(I_{H_1}+A_{H_!})
    \end{equation}
    If $A$ is trace class, the ordinary determinants also satisfy $ {\det}(I_H+A)={\det}(I_{H_1}+A_{H_!})$.
\end{lemma}
\begin{lemma}
    Let $\bm{a}\in\mathcal{A}_{N}$ and $g\in\Gamma_m$, and assume that $N_1\leq N_2\leq N,\ c_1\leq c_2\leq c$. Suppose that $(N_1,c_1)$ and $(N_2,c_2)$ satisfy \eqref{eq:condition of MNc} and that the two operators
    \begin{equation}
    \left\{
        \begin{aligned}
            T(\bm{a}):H_{N_1,c_1}\rightarrow H_{N_1,c_1}\\
            T(\bm{a}):H_{N_2,c_2}\rightarrow H_{N_2,c_2}
        \end{aligned}\right.
    \end{equation}
    are bijective. Then each of the determinants in \eqref{eq:def of tau} and \eqref{eq:def of tau2} has the same value on the two spaces. In particular, $\tau_{\bm{a}}(g),\tau^{(2)}_{\bm{a}}(g)$ are independent of $M,N,c$.
\end{lemma}
\begin{proof}
Write $H_i=H_{N_i,c_i}$ and let $T_i$ denote the realization of $T$ on $H_i$.
The weights give $H_1\subset H_2$, and $H(D_+)\subset H_1$ is dense
in $H_2$ by Lemma~2.1. By \eqref{eq:calculation of T(ga)-gT(a)}, the
Toeplitz operators agree on $H_1$. Bijectivity gives the same restriction
identity for their inverses. Hence
\begin{equation}
 \left.\bigl(g^{-1}T_2(g\bm{a})T_2(\bm{a})^{-1}-I\bigr)\right|_{H_1}
 =g^{-1}T_1(g\bm{a})T_1(\bm{a})^{-1}-I .
\end{equation}
The determinant identities now follow from Lemma~\ref{lem:subspace},
as in \cite[Lemma~15]{Kotani2023}.
\end{proof}
The class of symbols $\bm{a}$ with invertible Toeplitz operator is
\begin{equation}
    \mathcal{A}_{N}^{inv}=\{\bm{a}\in\mathcal{A}_{N}:T(\bm{a})^{-1}\text{ exist as a bounded operator on }H_{N,c}(D_+)\}
\end{equation}
Given $\bm{a}\in \mathcal{A}_{N},g\in\Gamma_m$, the product may satisfy $g\bm{a}\notin \mathcal{A}_{N}$. We therefore define the $\tau$-function for $g\bm{a}$ as follows.
\begin{lemma}
    Let $\bm{a}\in \mathcal{A}^{inv}_{N}$ and $g\in\Gamma_m$.\\
(i) $\tau_{\bm{a}}(g)$ is well defined, and $\tau_{\bm{a}}(g)\neq0$ if and only if $g^{-1}T(g\bm{a})$ is invertible on $H_{N,c}$ for every admissible $N$ and $c$.\\
(ii) $\tau_{g\bm{a}}^{(2)}(g_1)$ is given by
    \begin{equation}
        \tau_{g\bm{a}}^{(2)}(g_1)={\det}_2(g^{-1}g_1^{-1}T(g_1g\bm{a})(g^{-1}T(g\bm{a}))^{-1})
    \end{equation}
    and $\tau_{g\bm{a}}^{(2)}(g_1)\neq0$ if and only if $g^{-1}g_1^{-1}T(g_1g\bm{a})$ is invertible on $H_{N,c}(D_+)$ for every $c>0$.
\end{lemma}
To compute $\tau_{g\bm{a}}(g_1)$ when $g_1$ is rational, we first identify a finite-dimensional subspace containing the range of the perturbation of the identity in $g^{-1}g_1^{-1}T(g_1g\bm{a})(g^{-1}T(g\bm{a}))^{-1}$.
\begin{lemma}\label{finite rank}
    Suppose that $\tau_{\bm{a}}^{(2)}(g)\neq0$ for some $\bm{a}\in\mathcal{A}^{inv}_{N}$ and $g\in \Gamma_m$. For any rational function $r=\frac{p}{q}$ with neither poles nor zeros in $\overline{D}_+$, the operator $g^{-1}r^{-1}T(rg\bm{a})-g^{-1}T(g\bm{a})$ has finite rank, and its range $\operatorname{Ran}(g^{-1}r^{-1}T(rg\bm{a})-g^{-1}T(g\bm{a}))$ lies in the finite-dimensional space
    \begin{equation}
        span \{g (z)^{-1}r(z)^{-1}z^i,g (z)^{-1}z^j\}_{0\leq i\leq\deg p-1,0\leq j\leq\deg q-1}.
    \end{equation}
\end{lemma}

\begin{proof}
In \eqref{2.10} with multipliers $gr$ and $g$, the dependence of the
kernel on $z$ is given by
\[
 \frac{g(z)^{-1}}{q(\lambda)}
 \left(r(z)^{-1}\frac{p(\lambda)-p(z)}{\lambda-z}
             -\frac{q(\lambda)-q(z)}{\lambda-z}\right).
\]
The two divided differences are polynomials in $z$ of degrees at most
$\deg p-1$ and $\deg q-1$, respectively, and vanish for a constant
polynomial. Integrating in $\lambda$ gives the stated range inclusion
(cf.\ \cite[Lemma~37(i)]{Kotani2023}).
\end{proof}

\subsection{Properties of tau-functions}
Throughout this subsection, assume that
\begin{equation}
    \bm{a}\in \mathcal{A}^{inv}_{N},\ g_j=r_je^{h_j}\in\Gamma_m\text{with}\ \deg{r_j}=m_j\ \text{for}\ j=1,2.
\end{equation}
Set
\begin{equation}\label{Ea}
\begin{aligned}
    E_{\bm{a}}(g_1,g_2)&=\operatorname{tr}((g_1g_2)^{-1}T(g_1g_2\bm{a})T(g_1\bm{a})^{-1}g_1-I)(g_1^{-1}T(g_1\bm{a})T(\bm{a})^{-1}-I).\\
    &=\operatorname{tr}((g_1g_2)^{-1}R_{\bm{a}}(g_1,g_2)T(g_1\bm{a})^{-1}g_1)(g_1^{-1}R_{\bm{a}}(1,g_1)T(\bm{a})^{-1})\\
    &=\operatorname{tr}((g_1g_2)^{-1}R_{\bm{a}}(g_1,g_2)T(g_1\bm{a})^{-1}R_{\bm{a}}(1,g_1)T(\bm{a})^{-1})
\end{aligned}
\end{equation}
where $R_{\bm{a}}(g_1,g_2)=T(g_1g_2\bm{a})-g_2T(g_1\bm{a})$.

\begin{lemma}
    [Cocycle property] Let $g,g_1,g_2\in\Gamma_m$ and $\bm{a}\in\mathcal{A}_{N}$. If $\tau_{\bm{a}}^{(2)}(g_1)\neq0$, then
    \begin{equation}
        \tau_{\bm{a}}^{(2)}(g_1g_2)=\tau_{\bm{a}}^{(2)}(g_1)\tau_{g_1\bm{a}}^{(2)}(g_2)\exp{(-E_{\bm{a}}(g_1,g_2))}.
    \end{equation}
    Moreover, if $g_1=r_1$ and $g_2=r_2$ are rational functions, then 
    \begin{equation}\label{eq: r1r2}
        \tau_{g\bm{a}}(r_1r_2)=\tau_{g\bm{a}}(r_1)\tau_{r_1g\bm{a}}(r_2)
    \end{equation}
\end{lemma}

\begin{proof}
Factor the operator in \eqref{eq:def of tau2} as
\[
 (g_1g_2)^{-1}T(g_1g_2\bm{a})T(\bm{a})^{-1}
 =\bigl((g_1g_2)^{-1}T(g_1g_2\bm{a})T(g_1\bm{a})^{-1}g_1\bigr)
  \bigl(g_1^{-1}T(g_1\bm{a})T(\bm{a})^{-1}\bigr).
\]
Each factor differs from the identity by a Hilbert--Schmidt operator.
The product formula for $\det_2$ gives the exponential correction
$\exp(-E_{\bm{a}}(g_1,g_2))$ by \eqref{Ea}.
For rational factors, the same factorization at $g\bm{a}$ involves
finite-rank perturbations by Lemma~\ref{finite rank}. Multiplicativity
of the ordinary determinant gives \eqref{eq: r1r2}
(cf.\ \cite[Lemma~17(ii)]{Kotani2023}).
\end{proof}

\begin{lemma}
    If $g_1,g_2\in\Gamma_{m}$ and $g_1$ is even, that is, $g_1(z)=g_1(-z)$, then 
    \begin{equation}
        T(g_2g_1\bm{a})=T(g_2\bm{a})g_1,\qquad T(g_1\bm{a})=T(\bm{a})g_1.
    \end{equation}
    In particular, if $g_1=r_1$ and $g_2=r_2$ are rational functions, then
    \begin{equation}
        \tau_{r_1g\bm{a}}(r_2)=\tau_{g\bm{a}}(r_2)
    \end{equation}
\end{lemma}

\begin{proof}
Since $g_1$ is even, $(g_1u)_e=g_1u_e$ and $(g_1u)_o=g_1u_o$.
The two Toeplitz identities follow from the definition.
For $g_1=r_1$, they show that multiplication by $r_1$ intertwines the
finite-rank perturbations in the two determinants, as in
\cite[Lemma~17(iii)]{Kotani2023}. This multiplication is injective
and has dense range: its image is $H_{N+\deg r_1,c}(D_+)$,
which is dense by the approximation used in Lemma~2.1.
It therefore maps one finite-dimensional range bijectively onto the other.
Taking determinants on these ranges proves the assertion.
\end{proof}

We first compute the tau-function for rational factors, which we shall later use to approximate $e^h$. For $\zeta\in D_-$, let $q_\zeta(z)=\frac{1}{1-\zeta^{-1}z}$. This function is holomorphic in $D_+$.

\begin{lemma}
    Let $g\in\Gamma_m,\bm{a}\in\mathcal{A}_{N},\zeta,\zeta_1,\zeta_2\in D_-$ and suppose that $\tau_{\bm{a}}^{(2)}(g)\neq0$. Then
    \begin{equation}\left\{
    \begin{aligned}
        &\tau_{g\bm{a}}(q_\zeta)=1+\varphi_{g\bm{a}}^{(0)}(\zeta)\\
        &\tau_{g\bm{a}}(q_{\zeta_1}q_{\zeta_2})=\frac{(\zeta_1+\varphi^{(1)}_{g\bm{a}}(\zeta_1))(1+\varphi^{(0)}_{g\bm{a}}(\zeta_2))-(\zeta_2+\varphi^{(1)}_{g\bm{a}}(\zeta_2))(1+\varphi^{(0)}_{g\bm{a}}(\zeta_1))}{\zeta_1-\zeta_2}
    \end{aligned}\right.
    \end{equation}
\end{lemma}

\begin{proof}
By Lemma~\ref{finite rank}, for a single factor the perturbation has
range in the span of $g^{-1}$. For $\varphi\in H(D_-)$ we have
\begin{equation}\label{decomposition of qv}
        q_\zeta(z)\varphi(z)=q_\zeta(z)\varphi(\zeta)+q_\zeta(z)(\varphi(z)-\varphi(\zeta))\in H_{N,c}(D_+)\oplus H(D_-)
    \end{equation}
Thus the normalized operator acts on $g^{-1}$ by multiplication by
$1+\varphi^{(0)}_{g\bm{a}}(\zeta)$, giving the first identity.

For $q_{\zeta_1}q_{\zeta_2}$ the perturbation has range in the span
of $g^{-1}$ and $g^{-1}z$. For $j=1,2$ set
$a_j=1+\varphi^{(0)}_{g\bm{a}}(\zeta_j)$ and
$b_j=\zeta_j+\varphi^{(1)}_{g\bm{a}}(\zeta_j)$.
For $\zeta_1\ne\zeta_2$, subtract from each characteristic function
its linear interpolant at the two nodes. Its difference from the
interpolant, multiplied by $q_{\zeta_1}q_{\zeta_2}$, belongs to $H(D_-)$
because both poles cancel. Thus, after removing the common factor
$g^{-1}$, the images of $1,z$ have values $a_j,b_j$ at $\zeta_j$.
Taking determinants in evaluation coordinates gives
\[
 \tau_{g\bm{a}}(q_{\zeta_1}q_{\zeta_2})
 =\frac{\det\!\begin{pmatrix}a_1&b_1\\a_2&b_2\end{pmatrix}}
        {\det\!\begin{pmatrix}1&\zeta_1\\1&\zeta_2\end{pmatrix}}
 =\frac{a_2b_1-a_1b_2}{\zeta_1-\zeta_2},
\]
as in \cite[Lemma~37(ii)]{Kotani2023}. The coincident-node formula
follows by taking the limit.
\end{proof}

\subsection{Continuity of tau-functions}
Recall the metric \begin{equation}
    d_c(g_1,g_2)=\sqrt{\frac{1}{4\pi^2}\int_{C^2}\left|\frac{g_1^{-1}(z)g_1(\lambda)-g_2^{-1}(z)g_2(\lambda)}{\lambda-z}\right|^2|z^N\rho_c(z)|^{-2}|\lambda^N\rho_c(\lambda)|^{-2}|dz||d\lambda|}
\end{equation}
on $\Gamma_m$. The estimate $\|g_1^{-1}T(g_1\bm{a})-g_2^{-1}T(g_2\bm{a})\|_{H.S.}\leq c_{\bm{a}} d_c(g_1,g_2)$ gives the following continuity bound. 
\begin{lemma}\label{estimate for difference of tau function}
    Let $\bm{a}\in\mathcal{A}^{inv}_{N}$ and $g=re^h,g_1=r_1e^{h_1},g_2=r_2e^{h_2}\in\Gamma_m$, where $\deg h_j=m_j,\deg r_j=n_j$. Suppose that $M,N$ satisfy $M\geq m$ and $N>\max\{-n-1,n,\frac{m}{2}\}$, and assume $c>C$. Then there is a constant $c_{\bm{a}}$ such that
    \begin{equation}
        \left|\tau^{(2)}_{g\bm{a}}(g_1)-\tau^{(2)}_{g\bm{a}}(g_2)\right|\leq c_0d_{c}(gg_1,gg_2)
    \end{equation}
\end{lemma}
\begin{proof}
    If $\|A\|_{H.S.}, \|B\|_{H.S.}\leq c_1$, there is a constant $c_2$, depending only on $c_1$, such that 
    \begin{equation}
        \left|{\det}_{2}(I+A)-{\det}_2(I+B)\right|\leq c_2\|A-B\|_{H.S.}
    \end{equation}
    Since $\tau_{\bm{a}}^{(2)}(g)\neq0$, the operator $g^{-1}T(g\bm{a})$ has a bounded inverse. Applying the preceding estimate gives
    \begin{equation}\begin{aligned}
        \left|\tau^{(2)}_{g\bm{a}}(g_1)-\tau^{(2)}_{g\bm{a}}(g_2)\right|&\leq c_2\|(g^{-1}g_1^{-1}T(gg_1\bm{a})-g^{-1}g_2^{-1}T(gg_2\bm{a}))(g^{-1}T(g\bm{a}))^{-1}\|_{H.S.}\\
        &\leq \|(g^{-1}T(g\bm{a}))^{-1}\|c_{\bm{a}}c_2d_{c}(gg_1,gg_2)=c_0d_{c}(gg_1,gg_2)
        \end{aligned}
    \end{equation}
\end{proof}
\begin{lemma}
    Let $\bm{a}\in\mathcal{A}^{inv}_{N}$, with $M\geq m,N>n+1,c>c_2,\rho_c(z)=(\cosh z^M+K)^c$. Suppose that $g_0,g_k\in\Gamma_m$ with the rational function part of degree 0 satisfy
    \begin{equation}\label{in the sense of}
    \left\{
        \begin{aligned}
            &(i)\text{ there exist } c_1,c_2 \text{ such that} \\
            &\text{ for any }z\in C_x,|g^{-1}(z)g_k^{-1}(z)|,|g(z)g_k(z)|,|g(z)g_k'(z)|\leq c_1 \left|z^ne^{c_2|\Re{z}|^m}\right|\\
            &\text{ for any }z\in C_y,|g^{-1}(z)g_k^{-1}(z)|,|g(z)g_k(z)|,|g(z)g_k'(z)|\leq c_1 \left|z^n\right| \\
            &(ii)g_{k}(z)\rightarrow g_{0}(z)\text{ for any }z\in C 
        \end{aligned}\right.
    \end{equation}
    Then
    \begin{equation}
        \tau_{g\bm{a}}^{(2)}(g_k)\rightarrow\tau_{g\bm{a}}^{(2)}(g_0)
    \end{equation}
\end{lemma}
\begin{proof}
    By \cref{estimate for difference of tau function}, it suffices to prove that $d_{c}(gg_k,gg_0)\rightarrow0$. Set 
    \begin{equation}
        \Delta_k(z,\lambda)=\frac{g^{-1}(z)g(\lambda)g_k^{-1}(z)g_k(\lambda)-g^{-1}(z)g(\lambda)g_0^{-1}(z)g_0(\lambda)}{\lambda-z},
    \end{equation}
    \begin{equation}
        d_{c}^2(gg_k,gg_0)=\frac{1}{4\pi^2}\int_{C^2}|\Delta_k(z,\lambda)|^2|z^N\rho_c(z)|^{-2}|\lambda^N\rho_c(\lambda)|^{-2}|dz||d\lambda|,
    \end{equation}
    then we have
    \begin{equation}
    \begin{aligned}
        |\Delta_k(z,\lambda)|&\leq \left|g^{-1}(z)g(\lambda)g_k^{-1}(z)\frac{(g_k(\lambda)-g_k(z))}{\lambda-z}-g^{-1}(z)g(\lambda)g_{0}^{-1}(z)\frac{g_0(\lambda)-g_0(z)}{\lambda-z}\right|\\
        &\leq|g^{-1}(z)g^{-1}_{k}(z)||g(\lambda)\int_{0}^{1}g_k'((\lambda-z)t+z)dt|+|g^{-1}(z)g^{-1}_{0}(z)||g(\lambda)\int_{0}^{1}g_0'((\lambda-z)t+z)dt|\\
        &\leq 2c_3\left|\lambda^{n}z^ne^{c_2(|\Re{z}|^m+|\Re{\lambda}|^m)}\right|
    \end{aligned}
    \end{equation}
    on the region $|\lambda-z|\leq\varepsilon|\lambda|$ and
    \begin{equation}
        \begin{aligned}
            |\Delta_k(z,\lambda)|&= \left|\frac{g^{-1}(z)g(\lambda)g_k^{-1}(z)g_k(\lambda)-g^{-1}(z)g(\lambda)g_0^{-1}(z)g_0(\lambda)}{\lambda-z}\right|\\
            &\leq 2c_4\varepsilon^{-1}\left|\lambda^{n-1}z^ne^{c_2(|\Re{z}|^m+|\Re{\lambda}|^m)}\right|
        \end{aligned}
    \end{equation}
    on the region $|\lambda-z|\geq\varepsilon|\lambda|$. Together these estimates give 
    \begin{equation}
        |\Delta_k(z,\lambda)|^2\leq4\left|\lambda^{n}z^ne^{c_2(|\Re{z}|^m+|\Re{\lambda}|^m)}\right|^2(c_3^2\chi_{|\lambda-z|\leq\varepsilon|\lambda|}+c_4^2\varepsilon^{-2}|\lambda|^{-2}\chi_{|\lambda-z|>\varepsilon|\lambda|}).
    \end{equation}
    For $c_2<c$, the right-hand side is integrable on $(C^2,|\lambda^N\rho_c(\lambda)|^{-2}|z^N\rho_c(z)|^{-2}|dz||d\lambda|)$. Dominated convergence therefore gives 
    \begin{equation}
        d_c^2(gg_k,gg_0)\rightarrow\frac{1}{4\pi^2}\int_{C^2}\lim_{k\rightarrow\infty}|\Delta_k(z,\lambda)|^2|z^N\rho_c(z)|^{-2}|\lambda^N\rho_c(\lambda)|^{-2}|dz||d\lambda|=0
    \end{equation}

\end{proof}

\section{The nonnegativity condition defining
  \texorpdfstring{$\mathcal{A}_{N,+}^{inv}$}{A-inv(N,+)}}
We call $\bm{a}\in\mathcal{A}_{N},g\in\Gamma_m$ real if they satisfy
\begin{equation}
    \bm{a}(\bar{\lambda})=\overline{a(\lambda)}\ for\ \lambda\in C,\qquad\ \overline{g(z)}=g(\bar{z})\ for\ z\in\mathbb{C}
\end{equation}

Set 
\begin{equation}
\left\{
\begin{aligned}
    &\Gamma_{m}^{real}(C)=\{g\in\Gamma_{m}(C):g=\bar{g}\},\\
    &\mathcal{A}_{N,+}^{inv}=\{\bm{a}\in \mathcal{A}^{inv}_{N}:\bm{a}=\bar{\bm{a}}\text{ and }\tau_{\bm{a}}(r)\geq0\text{ for any real rational }\text{function }r\text{ with no poles nor zeros in }\overline{D}_{+}\}.
\end{aligned}\right.
\end{equation}

For rational factors, $\tau_{\bm{a}}(r)\geq0\iff\tau^{(2)}_{\bm{a}}(r)\geq0$ and $\tau_{\bm{a}}(r)>0\iff\tau^{(2)}_{\bm{a}}(r)>0$.

\begin{lemma}\label{rneh}
    If $\bm{a}\in \mathcal{A}_{N,+}^{inv}$ and $g\in\Gamma_{m}^{real}$, then $\tau_{\bm{a}}^{(2)}(g)\geq0$. Moreover, if $\tau_{\bm{a}}^{(2)}(g)>0$, then $\tau_{g\bm{a}}^{(2)}(r)\geq0$ for every rational function $r\in\Gamma_{0}^{real}$ with degree 0.
\end{lemma}
\begin{proof}
    Let $e^h\in\Gamma_{m}^{real}$, where $h(z)=\sum_{k=1}^{m}h_kz^k$ is an odd polynomial. We use the approximation
    
    \begin{equation}
       e^{h}=e^{h(z)-z^{4m}}e^{z^{4m}}=\lim_{n\rightarrow\infty}\left(1-\frac{h(z)-z^{4m}}{n}\right)^{-n}\left(1+\frac{z^{4m}}{n}\right)^n \text{ on }D_+
    \end{equation}
    and set $r_n=\left(1-\frac{h(z)-z^{4m}}{n}\right)^{-n}\left(1+\frac{z^{4m}}{n}\right)^n=r_{\zeta_1}\cdots r_{\zeta_{m_1}}r^{-1}_{\eta_1}\cdots r^{-1}_{\eta_{m_2}}$, where $\zeta_j,\eta_j\in D_-$.
    
    For sufficiently large $n$, the zeros of $r_n$ lie in $D_-$. To locate the poles, set $z=n^{\frac{1}{4m}}w$. Then $h(z)=O(n^{\frac{1}{4}}w^{m})$, and the rescaled pole equation is $w^{4m}+1+\frac{1}{n}h(n^{\frac{1}{4m}}w)=0$. Choose $R>0$ so that $|w|^{4m}-|h_1||w|^m>1$ for $|w|>R$. The zeros $\{w_{n,j}\}_{j=1}^{4m}$ of $w^{4m}+1+\frac{1}{n}h(n^{\frac{1}{4m}}w)$ lie in the disk $|w|\leq R$ and can be indexed so that $w_{n,j}\rightarrow w_j=e^{\frac{(2j+1)\pi i}{4m}}$ as $n\rightarrow\infty$. Hence, for sufficiently large $n$, $|\arg{z_{n,j}}-\frac{k\pi}{2}|>\frac{\pi}{8m}$. Thus $r_n$ has neither poles nor zeros in $D_+$.

    We next verify that $r_n\rightarrow e^{h}$ in the sense of \eqref{in the sense of}. For $|\Re{z}|>a>0$, we have 
    \begin{equation}\label{eq:estimate of rn x}\left\{
        \begin{aligned}
            &|r_n(z)|\leq\left(1+\left|\frac{h(z)}{n-h(z)+z^{4m}}\right|\right)^n\leq\left(1+\frac{|h(z)|}{n}\right)^n\leq e^{|h(z)|}\leq e^{c_2|\Re{z}|^m}\\
            &|r_n^{-1}(z)|\leq\left(1+\left|\frac{h(z)}{n+z^{4m}}\right|\right)^n\leq\left(1+\frac{|h(z)|}{n}\right)^n\leq e^{|h(z)|}\leq e^{c_2|\Re{z}|^m}\\
             &|r_n'(z)|=|r_n|\left|\frac{4mnz^{4m-1}}{n+z^{4m}}-\frac{n(4mz^{4m-1}-h'(z))}{n+z^{4m}-h(z)}\right|\leq 9m|z|^{4m-1}e^{|h(z)|}\leq e^{c_2|\Re{z}|^m}\\
             &|(r_n^{-1})'(z)|=|r_n^{-1}(z)|\left|\frac{n\left(4mz^{4m-1}-h'(z)\right)}{n+z^{4m}-h(z)}-\frac{4mnz^{4m-1}}{n+z^{4m}}\right|\leq9m|z|^{4m-1}e^{|h(z)|}\leq e^{c_2|\Re{z}|^m}
        \end{aligned}\right.
    \end{equation}
    whereas for $|\Re{z}|\leq a$ we have
       \begin{equation}\label{eq: estimate of rn y}\left\{
        \begin{aligned}
           &|r_{n}(z)|\leq\left|\frac{n+z^{4m}}{n+z^{4m}-h(z)}\right|^n\lesssim\left|\frac{(n+y^{4m})^2}{(n+y^{4m})^2+|h(iy)|^2}\right|^{\frac{n}{2}}\leq c_1\\
           &|r_{n}^{-1}(z)|\leq\left|\frac{n+z^{4m}-h(z)}{n+z^{4m}}\right|^n\lesssim\left|\frac{(n+y^{4m})^2+|h(iy)|^2}{(n+y^{4m})^2}\right|^{\frac{n}{2}}\leq \left|1+\frac{|2h_1y^m|^2}{2ny^{4m}}\right|^{\frac{n}{2}}
        \leq e^{h_1^2y^{-2m}}\leq c_1\\
        &|r_n'(z)|=|r_n|\left|\frac{4mnz^{4m-1}}{n+z^{4m}}-\frac{n(4mz^{4m-1}-h'(z))}{n+z^{4m}-h(z)}\right|\leq c_1|z|^{4m-1}\\
        &|(r_n^{-1})'(z)|=|r_n^{-1}(z)|\left|\frac{n\left(4mz^{4m-1}-h'(z)\right)}{n+z^{4m}-h(z)}-\frac{4mnz^{4m-1}}{n+z^{4m}}\right|\leq c_1|z|^{4m-1}
        \end{aligned}\right.
    \end{equation}
    It follows that $\tau_{\bm{a}}^{(2)}(g)=\lim_{n\rightarrow\infty}\tau_{\bm{a}}^{(2)}(r_n)\geq0$.

    If $\tau_{\bm{a}}^{(2)}(g)>0$, the cocycle identity gives
    \begin{equation}
        \tau_{g\bm{a}}^{(2)}(r)=\frac{\tau_{a}^{(2)}(gr)}{\tau_{a}^{(2)}(g)}\exp(E_{\bm{a}}(g,r))\geq0
    \end{equation}
    for every $r\in\Gamma_{0}^{real}$ with $\deg r=0$.
\end{proof}

\begin{lemma}\label{lem:tau r>0}
    Let $\bm{a}\in\mathcal{A}^{inv}_{N,+}$ and $g\in\Gamma_{m}^{real}$, with $M\geq m$, and suppose that $\tau_{\bm{a}}(g)>0$. Then the following assertions hold.\\
(i) The function $m_{g\bm{a}}$ satisfies
    \begin{equation}
        \frac{\Im{m_{g\bm{a}}}}{\Im{z}}>0
    \end{equation}
    and $1+\varphi^{(0)}_{g\bm{a}}$ has no zeros in $D_-$.\\
(ii) For every $r\in\Gamma_0^{real}$ with $\deg r=0$, we have $\tau_{g\bm{a}}^{(2)}(r)>0$.
\end{lemma}
\begin{proof}
    Recall that $\tau_{g\bm{a}}(q_{\zeta_1}q_{\zeta_2})=\frac{m_{g\bm{a}}(\zeta_1)-m_{g\bm{a}}(\zeta_2)}{\zeta_1-\zeta_2}\left(1+\varphi_{g\bm{a}}^{(0)}(\zeta_1)\right)\left(1+\varphi_{g\bm{a}}^{(0)}(\zeta_2)\right)$ for all $\zeta_1,\zeta_2\in D_-$. For a fixed $\zeta\in D_-$, this gives 
    \begin{equation}
        0\leq\tau_{g\bm{a}}(q_{\zeta}q_{\bar{\zeta}})=\frac{\Im{m_{g\bm{a}}(\zeta)}}{\Im{\zeta}}\left|1+\varphi_{g\bm{a}}^{(0)}(\zeta)\right|^2.
    \end{equation}
    
    Set 
    \begin{equation}
        \mathcal{Z}=\left\{\zeta\in D_-:1+\varphi_{g\bm{a}}^{(0)}(\zeta)=0\right\}.
    \end{equation}
    Since $1+\varphi_{g\bm{a}}^{(0)}(\zeta)\rightarrow1$ as $\zeta\rightarrow\infty$, the function $1+\varphi_{g\bm{a}}^{(0)}(\zeta)$ is not identically $0$. Thus $1+\varphi_{g\bm{a}}^{(0)}(\zeta)$ is meromorphic in $D_-$, and $\frac{\Im{m_{g\bm{a}}(\zeta)}}{\Im{\zeta}}\geq0$ on $D_-\setminus\mathcal{Z}$. Suppose that $\zeta_0\in D_-\cap\mathcal{Z}\cap \mathbb{C}_+$. Then $\Im{\left(-m_{g\bm{a}}(\zeta_0)^{-1}\right)}=0$, so the minimum principle gives $\Im{\left(-m_{g\bm{a}}(\zeta)^{-1}\right)}=0$ throughout the component of $D_-\cap \mathbb{C}_+$ containing $\zeta_0$. This contradicts $m_{g\bm{a}}(\zeta)-\zeta\rightarrow0$ as $\zeta\rightarrow\infty$. Thus $D_-\cap\mathcal{Z}=\varnothing$ and $\frac{\Im{m_{g\bm{a}}(\zeta)}}{\Im{\zeta}}\geq0$ on $D_-$. Applying the minimum principle once more, this time to $\Im{m_{g\bm{a}}(z)}$, gives $\frac{\Im{m_{g\bm{a}}}}{\Im{z}}>0$.

    For (ii), we follow the proof of \cite[Lemma~4.16]{Kotani2024}. Since $D_-\cap \mathbb{R}=\varnothing$, it suffices to consider $r=q_{\zeta}q_{\bar{\zeta}}(q_{\eta}q_{\bar{\eta}})^{-1}$. The argument applies on each of the two connected components of $\mathbb{C}_+\cap D_-$.
\end{proof}

\begin{lemma}\label{lem:tau ga(g)>0}
    Let $\bm{a}\in\mathcal{A}^{inv}_{N,+}$ and $g\in\Gamma_{m}^{real}$, with $M\geq m$, and suppose that $\tau_{\bm{a}}^{(2)}(g)>0$. Then, for every $re^h\in\Gamma_{m}^{real}$ with $\deg r=0$, we have $\tau_{g\bm{a}}^{(2)}(re^h)>0$.
\end{lemma}
\begin{proof}
    Choose rational approximants with $r_k\rightarrow re^h$ in the sense of \eqref{in the sense of}. The calculation in \cref{rneh} gives $r_k^{-1}re^h\rightarrow1$. Hence $\tau_{g\bm{a}}^{(2)}(re^hr_k^{-1})>0$ for all sufficiently large $k$. Lemma~\ref{lem:tau r>0} then gives $\tau_{re^h r_k^{-1}g\bm{a}}^{(2)}(r_k)>0$, and the cocycle identity yields
    \begin{equation}
        \tau_{g\bm{a}}^{(2)}(re^h)=\tau_{g\bm{a}}^{(2)}(re^hr_k^{-1}r_k)=\tau_{g\bm{a}}^{(2)}(re^hr_k^{-1})\tau_{re^hr_k^{-1}g\bm{a}}^{(2)}(r_k)\exp(-E_{g\bm{a}}(re^hr_k^{-1},r_k))>0
    \end{equation}
\end{proof}

\section{Derivation of the Schr\"odinger and KdV equations}
We derive the Schr\"odinger equation and the KdV evolution from the group elements $e_x(z)=e^{xz}, e_{t,x}=e^{tz^3+xz}$.

Assume that $e_{t,x}\bm{a}\in\mathcal{A}_{N}^{inv}$ for all $x,t\in\mathbb{R}$, and define 
\begin{equation}
    u_{x,t}=T(e_{x,t}\bm{a})^{-1}1\in H_{1}(D_+),\qquad w_{x,t}=\mathfrak{p}_-(e_{x,t}\bm{a}u_{x,t})\in H_{-}.
\end{equation}
Choose a bounded analytic function $\tilde{\bm{a}}(z)$ on $D_+$ such that $\rho_c(\lambda)(\tilde{\bm{a}}(\lambda)-\bm{a}(\lambda))$ is bounded on the contour. At $t=0$, the Cauchy representation is 
\begin{equation}
    w_x(z)=\frac{1}{2\pi i}\int_C\frac{e^{x\lambda}(\tilde{\bm{a}}(\lambda)-\bm{a}(\lambda))u_{x}(\lambda)}{z-\lambda}d\lambda
\end{equation}
Differentiation with respect to $x$ gives
\begin{equation}\begin{aligned}
    &\partial_x w_x(z)=\frac{1}{2\pi i}\int_{C}e^{x\lambda}(\bm{a}(\lambda)-\tilde{\bm{a}}(\lambda))\frac{\lambda u_x(\lambda)+\partial_xu_{x}(\lambda)}{z-\lambda}\dd\lambda\\
    &\partial_x^2w_x(z)=\frac{1}{2\pi i}\int_{C}e^{x\lambda}(\bm{a}(\lambda)-\tilde{\bm{a}}(\lambda))\frac{\partial_x^2u_{x}(\lambda)+2\lambda\partial_x u_x(\lambda
    )+\lambda^2u_x(\lambda)}{z-\lambda}\dd\lambda
\end{aligned}
\end{equation}
since $\partial_x^ju_x\in H_{j+1}(D_+)$ for $j\leq N$ by \cite[Lemma~4.7]{Kotani2024}.

Since $u_x\in H_{1}(D_+)$, the expansion
\begin{equation}
    (z-\lambda)^{-1}=\sum_{1\leq k\leq L}z^{-k}\lambda^{k-1}+z^{-L}\lambda^{L}(z-\lambda)^{-1}
\end{equation}
and the boundedness of $|z^N\rho_{c}(z)||\tilde{\bm{a}}(z)-\bm{a}(z)|$ on $C$ give
\begin{equation}\begin{aligned}
    w_{x}(z)&=\sum_{1\leq k\leq N-1}z^{-k}\frac{1}{2\pi i}\int_{C}\lambda^{k-1}e^{x\lambda}(\bm{a}(\lambda)-\tilde{\bm{a}}(\lambda))\dd\lambda+z^{-N+1}\frac{1}{2\pi i}\int_C\frac{\lambda^{N-1}e^{x\lambda}(\tilde{\bm{a}}(\lambda)-\bm{a}(\lambda))u_{x}(\lambda)}{z-\lambda}d\lambda\\
    &=\sum_{1\leq k\leq N-1}z^{-k}s_k(x)+z^{-N+1}\tilde{w}_x(z)
    \end{aligned}
\end{equation}
with $\tilde{w}_x\in H_-$. For the first two derivatives, this gives
\begin{equation}
    \begin{aligned}
        \partial_x w_x(z)=\sum_{1\leq k\leq N-2}z^{-k}s_k(x)+z^{-N+2}\tilde{w}_x^{(1)}(z)\in H_-,\\
        \partial_x^2 w_x(z)=\sum_{1\leq k\leq N-3}z^{-k}s_k(x)+z^{-N+3}\tilde{w}_x^{(1)}(z)\in H_-.
    \end{aligned}
\end{equation}
The same calculation applies with $e_{x,t}$ in place of $e_x$.

In particular,
\begin{equation}
    \varphi_{e_{x,t}}(z)=w_{x,t}(z)=\sum_{k=1}^{\infty}z^{-k}s_k(x,t),s_1(x,t)=\kappa_{1}(e_{x,t}\bm{a}).
\end{equation}
For $x,t\in\mathbb{R}$, set 
\begin{equation}
    f_{\bm{a}}(x,t,z)=\bm{a}(z)u_{x,t}(z)=e_{x,t}(z)^{-1}(1+\varphi_{e_{x,t}}(z))
\end{equation}
Direct differentiation gives
\begin{equation}
    e_{x,t}(-\partial_x^2f_{\bm{a}}+qf_{\bm{a}}+z^2f_{\bm{a}})=-\partial_x^2\varphi_{e_{x,t}}+q(1+\varphi_{e_{x,t}})+2z\partial_{x}\varphi_{e_{x,t}}.
\end{equation}
The constant term is $q+2s_1'$. Choosing $q(x)=-2s_1'(x)$ gives 
\begin{equation}
    e_{x,t}(-\partial_x^2f_{\bm{a}}+qf_{\bm{a}}+z^2f_{\bm{a}})\in H_-
\end{equation}
The left-hand side also has the representation
\begin{equation}
    e_{x,t}\bm{a}(-\partial_x^2u_{x,t}+qu_{x,t}+z^2u_{x,t})\in (e_{x,t}\bm{a})H_{3}(D_+)
\end{equation}
since $u_{x,t},\partial_x^2 u_{x,t}\in H_{1}(D_+)$. By the assumed invertibility of $T(e_{x,t}\bm{a})$ on $H_{N}(D_+)$, we obtain 
\begin{equation}
    -\partial_x^2u_{x,t}+qu_{x,t}+z^2u_{x,t}=0.
\end{equation}
Multiplication by $\bm{a}$ gives the Schr\"odinger equation 
\begin{equation}
    -\partial_x^2f_{\bm{a}}(x,t,z)+q(x)f_{\bm{a}}(x,t,z)=-z^2f_{\bm{a}}(x,t,z).
\end{equation}
The Schr\"odinger equation gives 
\begin{equation}\begin{aligned}
        0&=-\partial_x^2\varphi_{e_{x,t}}+q(1+\varphi_{e_{x,t}})+2z\partial_{x}\varphi_{e_{x,t}}\\
        &=\sum_{1\leq k\leq N-3}(-s_k''-2s_1's_k+2s_{k+1}')z^{-k}+z^{-N+3}\tilde{w}_x,
\end{aligned}
\end{equation}
Comparing coefficients, we obtain the vanishing of 
\begin{equation}\label{eq:relations between sk}
    -s_k''-2s_1's_k+2s_{k+1}'\text{ for }1\leq k\leq N-3
\end{equation}
The coefficients are thus determined recursively by $s_1$ and the initial values $\{s_{k}(0)\}_{2\leq k\leq N-2}$.

We now determine the time evolution of the reconstructed potential.
\begin{lemma}\label{lem:KdV derivative}
    The potential $q(x,t)=-2s_1'(x,t)$ satisfies 
    \begin{equation}
        \partial_tq(x,t)=\frac{1}{4}\partial_x^3q(x,t)-\frac{3}{2}q(x,t)\partial_xq(x,t).
    \end{equation}
\end{lemma}
\begin{proof}
As in \cite[Proposition~14]{Kotani2023}, take
\begin{equation}
    L=-\partial_x^2+q,\qquad P=\partial_x^3-\frac{3}{2}q\partial_x-\frac{3}{4}q'.
\end{equation}
Since $Lf_{\bm{a}}=-z^2f_{\bm{a}}$, we have
$Pf_{\bm{a}}=(z^2-q/2)\partial_x f_{\bm{a}}+q'f_{\bm{a}}/4$.
Writing $e_{x,t}f_{\bm{a}}=1+\varphi$ as above gives
\[
 e_{x,t}(\partial_t f_{\bm{a}}-Pf_{\bm{a}})
 =\partial_t\varphi-(z^2-q/2)\varphi'
  -\tfrac12 zq(1+\varphi)-\tfrac14 q'(1+\varphi).
\]
The coefficients of $z$ and $1$ are $-s_1'-q/2$ and
$-s_2'-qs_1/2-q'/4$. Both vanish by $q=-2s_1'$ and
\eqref{eq:relations between sk} with $k=1$, so the residual belongs
to $H_-$. It also equals
$e_{x,t}\bm{a}(\partial_tu_{x,t}-(z^2-q/2)\partial_xu_{x,t}-q'u_{x,t}/4)$.
Projection and injectivity of the Toeplitz operator give
\begin{equation}
    \partial_t f_{\bm{a}}-Pf_{\bm{a}}=0 .
\end{equation}
Compatibility with $Lf_{\bm{a}}=-z^2f_{\bm{a}}$ gives
$(\partial_t L-[P,L])f_{\bm{a}}=0$. The commutator is multiplication
by $\tfrac14\partial_x^3 q-\tfrac32 q\partial_x q$. At each $(x,t)$
we may take $z$ large in an exterior sector so that
$1+\varphi(z)\ne0$. The asserted equation follows.

\end{proof}

\section{Proofs of the main results on the flow}

First we give the proof of Theorem \ref{Thm: theorem1}.

\begin{proof}[Proof of Theorem 1]
    For a given symbol $\bm{a}\in\mathcal{A}_{N,+}^{inv}$, we have defined the Toeplitz operator $T(g\bm{a})$ for any $g\in\Gamma_{m}$. Combining with Lemma \ref{lem:tau ga(g)>0} we obtain that $g^{-1}T(g\bm{a})$ is invertible on $H_{N,+}(D_+)$ for any $g=re^h\in \Gamma_{m}^{real}$ with $\deg r=0$. Thus, we can define the characteristic function of $e_{t,x}\bm{a}$ for any $x,t\in\mathbb{R}$.  The calculation in \cite[Lemma 4.7]{Kotani2024} yields the continuity of $q(x,t)$ both in $x$ and $t$. Finally applying Lemma \ref{lem:KdV derivative} we can obtain the solution for KdV equation. 
    
\end{proof}

Let $q$ be a real potential with half-line Weyl functions $m_\pm$. We give conditions on $m_{\bm{a}}$ under which $(1,\frac{m_a(z)}{z})\in\mathcal{A}_{N,+}^{inv}$, where
    
\begin{equation}
    m_{\bm{a}}(\zeta)=\left\{\begin{aligned}
        -m_+(-\zeta^2)\ \ \  \text{if}\ \Re{\zeta}>0\\
        m_-(-\zeta^2)\ \ \ \text{if}\ \Re{\zeta}<0
    \end{aligned}\right.
\end{equation}

\subsection{Estimates for spectral measures}
For $p>0$, define the following class of measures, with the integrability condition imposed for every $c>0$.
\begin{equation}
    \mathcal{E}_p=\left\{\nu:\int_{-\infty}^{0}e^{c|\lambda|^{p}}\nu(\dd\lambda)<\infty\right\}.
\end{equation}

    We first characterize membership in this class by a tail estimate. Suppose that $\nu\in\mathcal{E}_p$. 

    Given $a>0$, take $c=a$ in the moment condition. Then
    \begin{equation}\label{eq:tail estimate}
        \nu((-\infty,-R])\leq \int_{-\infty}^{-R}e^{-aR^p}e^{c|\lambda|^p}\nu(d\lambda)\leq C_ae^{-aR^p}\text{ for }R\geq1
    \end{equation}
    Conversely, suppose that \eqref{eq:tail estimate} holds for every $a>0$. To prove $\nu\in\mathcal{E}_p$, proceed as follows.  
    
    Set $R_n=2^nR_0$. For fixed $c>0$, choose $a>2^pc$. On each dyadic interval we have
    \begin{equation}
        \int_{-2R_n}^{-R_n}e^{c|\lambda|^p}\nu(\dd\lambda)\leq C_ae^{-(a-2^pc)R_n^p}
    \end{equation}
    Since $\sum_{n=1}^{\infty}C_a^{-(a-2^pc)(2^nR_0)^p}<\infty$ and $\nu([-R_0,0])<\infty$, summing over the intervals gives $\nu\in\mathcal{E}_p$. More generally, suppose that there is a measure $\eta\in\mathcal{E}_p$ such that, for all sufficiently large $R$,
    \begin{equation}
        \nu([-2R,-R])<CR^3\sqrt{\eta(-\infty,-R/4)}
    \end{equation}
    Then $\nu\in \mathcal{E}_p$, since the right-hand side is bounded by $CC_a^{1/2}R^3e^{-aR^p/2\cdot 4^p}$.

    Recall the integral representation of a Herglotz function $F$,
    \begin{equation}\label{eq:Herglotz expression of F}
        F(z)=\alpha+\beta z+\int_{\mathbb{R}}\left(\frac{1}{\lambda-z}-\frac{\lambda}{1+\lambda^2}\right)\nu_F(\dd\lambda)
    \end{equation}
    where $\int_{\mathbb{R}}(1+\lambda^2)^{-1}\nu_F(\dd\lambda)<\infty$.
    
    Write 
    \begin{equation}
        I_R=[-2R,-R],\qquad J_R=[-3R,-R/2].
    \end{equation}
    For $R\geq2,0\leq y\leq1$, the Poisson kernel gives the estimate
    \begin{equation}\label{eq:Poisson kernel}
        \nu_F(I_R)\leq\frac{2}{\pi}\int_{J_R}\Im F(x+iy)\dd x.
    \end{equation}
    Indeed, for each $\lambda\in[-2R,R]$,
    \begin{equation}
        [\lambda-R/2,\lambda+R/2]\subset[-3R,-R/2]=J_R
    \end{equation}
    and hence
    \begin{equation}
        \int_{J_R}\frac{y}{(\lambda-x)^2+y^2}\dd x\geq\int_{\lambda-R/2}^{\lambda+R/2}\frac{y}{(\lambda-x)^2+y^2}\dd x=2\arctan\frac{R}{2y}\geq\frac{\pi}{2}
    \end{equation}
    since $R\geq2\geq 2y$.
    
    Taking imaginary parts in \eqref{eq:Herglotz expression of F} gives 
    \begin{equation}
        \Im F(x+iy)=\beta y+\int_{\mathbb{R}}\frac{y\dd x}{(\lambda-x)^2+y^2}
    \end{equation}
    Dropping the nonnegative term $\beta y$ and restricting the integral to $I_R$, we obtain
    \begin{equation}\begin{aligned}
        \int_{J_R}\Im F(x+iy)\dd x&\geq\int_{J_R}\int_{I_R}\frac{y}{(\lambda-x)^2+y^2}\nu_F(\dd\lambda)\dd x\\
        &\geq\int_{I_R}\left(\int_{J_R}\frac{y}{(\lambda-x)^2+y^2}\dd x\right)\nu_F(\dd\lambda)\geq\frac{\pi}{2}\nu_F(I_R)
    \end{aligned}
    \end{equation}
    by Tonelli's theorem. We use this estimate to prove the following equivalence.
\begin{lemma}\label{comparison of nuF nu-1/F}
If $F$ is a Herglotz function, then $-1/F$ is also a Herglotz function, and for every $p>0$ we have
    \begin{equation}
        \nu_F\in\mathcal{E}_p\iff\nu_{-1/F}\in\mathcal{E}_p.
    \end{equation}
\end{lemma}
\begin{proof}

    Since $F$ is a Herglotz function, $F\neq0$ in $\mathbb{C}_+$. The identity
    \begin{equation}
        \Im\left(\frac{-1}{F}\right)=\frac{\Im F}{|F|^2}>0
    \end{equation}
    shows that $-1/F$ is also a Herglotz function. Suppose that $\nu:=\nu_F\in\mathcal{E}_p$. In particular, $\nu((-\infty,0))<\infty$, so we may rewrite \eqref{eq:Herglotz expression of F} as 
    \begin{equation}
        F(z)=H(z)+\int_{\mathbb{R}_-}\frac{\nu(\dd\lambda)}{\lambda-z}
    \end{equation}
    with
    \begin{equation}
        H(z)=\alpha'+\beta z+\int_{\mathbb{R}_+}\left(\frac{1}{\lambda-z}-\frac{\lambda}{1+\lambda^2}\right)\nu(\dd\lambda)
    \end{equation}
    where $\alpha'=\alpha+\int_{\mathbb{R}_-}\frac{-\lambda}{1+\lambda^2}\nu(\dd\lambda)$. The function $H$ is analytic across $\mathbb{R}_-$, is real there, and satisfies 
    \begin{equation}
        H'(x)=\beta+\int_{\mathbb{R}_+}\frac{\nu(\dd\lambda)}{(\lambda-x)^2}\geq0
    \end{equation}
    Hence the limit $l:=\lim_{r\rightarrow\infty}H(-r)$ exists in $[-\infty,\infty)$. Define
    \begin{equation}
        \delta_R=\nu((-\infty,-R/4)),\qquad F_R(z)=H(z)+\int_{-R/4}^{0}\frac{\nu(\dd\lambda)}{\lambda-z}.
    \end{equation}
    
    We divide the proof into three steps.

    \emph{Step 1.} For all sufficiently large $R$, 
    \begin{equation}\label{eq:lower bound of F_R(x)}
        |F_R(x)|\geq\frac{c_0}{R},\qquad x\in J_R
    \end{equation} with a constant $c_0>0$ independent of $R$.
    
    If $l\neq0$, there is an $R_0>0$ such that, for every $R>R_0$, we have $|H(x)|>\varepsilon$ throughout $x\in J_R=[-3R,-R/2]$. Since $\lambda-x\geq R/4$ on the integration interval, we also have 
    \begin{equation}\label{eq:4/R}
        0\leq\int_{-R/4}^{0}\frac{\nu(\dd\lambda)}{\lambda-z}\leq\frac{4}{R}\nu((-\infty,0))\rightarrow0.
    \end{equation}
    Hence $|F_R(x)|$ is bounded away from zero uniformly on $J_R$ for all sufficiently large $R$.

    If $l=0$, we start from the identity
    \begin{equation}
        H(-1)-H(-r)=\beta(r-1)+\int_{\mathbb{R}_+}\left(
        \frac{1}{\lambda+1}-\frac{1}{\lambda+r}\right)\nu(\dd\lambda).
    \end{equation}
    As $r\rightarrow+\infty$, the linear term would satisfy $|\beta(r-1)|\rightarrow\infty$ if $\beta>0$. Both terms on the right are nonnegative, and in particular 
    \begin{equation}
        \beta(r-1)+\int_{\mathbb{R}_+}\left(
        \frac{1}{\lambda+1}-\frac{1}{\lambda+r}\right)\nu(\dd\lambda)\geq0.
    \end{equation}
    
    Since $H(-1)-H(-r)$ has a finite limit as $r\rightarrow+\infty$, it follows that $\beta=0$. Monotone convergence then gives $\int_{\mathbb{R}_+}(1+\lambda)^{-1}\nu(\dd\lambda)<\infty$, and the assumption $l=0$ gives $\alpha'-\int_{\mathbb{R}_+}\frac{\lambda}{1+\lambda^2}\nu(\dd\lambda)=0$. Thus
    \begin{equation}
        H(z)=\int_{0}^{\infty}\frac{\nu(\dd\lambda)}{\lambda-z},\qquad F_R(z)=\int_{-R/4}^{\infty}\frac{\nu(\dd\lambda)}{\lambda-z}.
    \end{equation}
    In this case $\nu\neq0$, since otherwise the representation above would give $F=0$, contrary to the Herglotz property of $F$. Choose $K\geq1$ so that $\gamma:=\nu([-K,K])>0$. If $R\geq8K,x\in J_R$, then $0<\lambda-x<3R+K\leq25R/8$ for every $\lambda\in[-K,K]$. Consequently,
    \begin{equation}
        F_R(x)\geq\int_{_K}^{K}\frac{\nu(\dd\lambda)}{\lambda-x}\geq\frac{8\gamma}{25R}.
    \end{equation}
    \emph{Step 2.} For all sufficiently large $R$ and all $x\in J_R,0\leq y\leq1$, we have
    \begin{equation}\begin{aligned}
        |F_R'(x+iy)|&\leq\beta+\int_{\mathbb{R}_+}\frac{\nu(\dd\lambda)}{\abs{\lambda-x-yi}^2}+\int_{-R/4}^{0}\frac{\nu(\dd\lambda)}{\abs{\lambda-x-yi}^2}\\
        &\leq\beta+\int_{\mathbb{R}_+}\frac{\nu(\dd\lambda)}{(1+\lambda)^2}+\frac{16\nu((-\infty,0))}{R^2}\leq C.
    \end{aligned}
    \end{equation}
    The last term is estimated as in \eqref{eq:4/R}.

    For $x\in\mathbb{R}$ in $J_R$, the value $F_R(x)$ is real. Integrating the derivative along the vertical segment from $x$ to $x+iy$ gives
    \begin{equation}
        |F_R(x+iy)-F_R(x)|\leq Cy,\qquad 0\leq\Im F_R(x+iy)\leq Cy
    \end{equation}
    The omitted tail satisfies
    \begin{equation}\label{eq:estimate of FR,F}
        |F_R(x+iy)-F(x+iy)|=\abs{\int_{-\infty}^{-R/4}\frac{\nu(\dd\lambda)}{\lambda-x-yi}}\leq\frac{\delta_R}{y}
    \end{equation}
    Suppose that $\delta_R>0$ and set $y_R=\sqrt{\delta_R}$. By \eqref{eq:tail estimate}, we have $y_R\rightarrow0,Ry_R\rightarrow0$. Combining these bounds with \eqref{eq:lower bound of F_R(x)} gives
    \begin{equation}\begin{aligned}
                \abs{F(x+iy_R)}&\geq\abs{F_R(x)}-\abs{F_R(x+iy_R)-F_R(x)}-\abs{F_R(x+iy_R)-F(x+iy_R)}\\
                &\geq\frac{c_0}{R}-Cy_R-\frac{\delta_R}{y_R}\geq\frac{c_0}{2R}
    \end{aligned}
    \end{equation}
    while \eqref{eq:estimate of FR,F} gives
    \begin{equation}
        \Im F_R(x+iy_R)\leq C\sqrt{\delta_R},\qquad \Im F(x+iy_R)\leq C'\sqrt{\delta_R}
    \end{equation}
    Consequently, 
    \begin{equation}
        \Im\left(\frac{-1}{F_R(x+iy_R)}\right)\leq C''R^2\sqrt{\delta_R},
    \end{equation}
    Since $\abs{J_R}=5R/2$, applying \eqref{eq:Poisson kernel} to $-1/F$ gives
    \begin{equation}
        \nu_{-1/F}(I_R)\leq\frac{2}{\pi}\int_{J_R}\Im\left(\frac{-1}{F(x+iy_R)}\right)\dd x\leq \frac{5R}{2}\frac{2}{\pi}C''R^2\sqrt{\delta_R}\leq C'''R^3\sqrt{\nu_F((-\infty,-R/4))}.
    \end{equation}
    This proves $\nu_{-1/F}\in\mathcal{E}_p$.

    \emph{Step 3.} Conversely, suppose that $\nu_{-1/F}\in\mathcal{E}_p$. Applying Steps~1 and~2 with representing measure $\nu_{-1/F}$ gives the result, since $\nu_{F}=\nu_{-1/(-1/F)}$.
\end{proof}
Set $S=m_++m_-$, so that $\Im S>0$, and define
\begin{equation}\label{eq:def of V}
    \mathcal{V}=\frac{1}{S}\left(\begin{array}{cc}
         m_+m_-&(m_+-m_-)/2  \\
         (m_+-m_-)/2&-1 
    \end{array}\right)
\end{equation}
For $w\in\mathbb{C}$, write $e(w)=(w,1)^T$ and set $w_{\pm}=\Im m_{\pm},u=-m_+,v=m_-$. Then 
\begin{equation}
    \mathcal{V}_q=\frac{e(u)e(v)^T+e(v)e(u)^T}{2(u-v)}
\end{equation}
and
\begin{equation}
    \Im \mathcal{V}_q=\frac{w_+}{|S|^2}\left(\begin{array}{c}
         m_-  \\
          1
    \end{array}\right)\left(\begin{array}{c}
         m_-  \\
          1
    \end{array}\right)^*+\frac{w_+}{|S|^2}\left(\begin{array}{c}
         -m_+  \\
          1
    \end{array}\right)\left(\begin{array}{c}
         -m_+  \\
          1
    \end{array}\right)^*\geq0
\end{equation}
where $\Im\mathcal{V}_q=\frac{\mathcal{V}-\mathcal{V}^*}{2i}$. Hence $\mathcal{V}$ is a matrix-valued Herglotz function. Denote its representing measure by $\Sigma$ and set $\rho=\operatorname{tr}\Sigma=\Sigma_{11}+\Sigma_{22}$.
\begin{lemma}\label{lem:comparison of sigma and rho}
    The spectral measures $\sigma_\pm(\dd\lambda)$ of the half-line Weyl functions and the matrix spectral measure $\Sigma(\dd\lambda)$ of the Green function satisfy
    \begin{equation}
        \sigma_++\sigma_-\in\mathcal{E}_p\iff \rho=\operatorname{tr}\Sigma\in\mathcal{E}_p.
    \end{equation}
\end{lemma}

\begin{proof}
    Suppose first that $\sigma_++\sigma_-\in\mathcal{E}_p$. By Lemma~\ref{comparison of nuF nu-1/F}, $\Sigma_{22}=\nu_{-1/(m_++m_-)}\in\mathcal{E}_p$. For the first diagonal entry, $\Sigma_{11}=\frac{-1}{\frac{-1}{m_+}+\frac{-1}{m_-}}$ and $\nu_{\frac{-1}{m_+}+\frac{-1}{m_-}}\in\mathcal{E}_p$ give $\Sigma_{11}\in\mathcal{E}_p$ by the same lemma. Conversely, suppose that $\Sigma_{11}+\Sigma_{22}\in\mathcal{E}_p$. Since $\Sigma_{22}=\nu_{-1/(m_++m_-)}$, the lemma gives $\sigma_++\sigma_-=\nu_{m_++m_-}\in\mathcal{E}_p$.
\end{proof}

\subsection{Approximation of the initial data on the contour}
To distinguish the behavior near the real and imaginary axes, write $C$ as the union of the four arcs
\begin{equation}
    C_1=C\cap\{\Re{z}\geq x_0\},\qquad C_2=C\cap\{\Im{z}>y_0\},\qquad C_3=C\cap\{\Re{z}\leq-x_0\},\qquad C_4=C\cap\{\Im{z}<-y_0\}, 
\end{equation}
where $x_0,y_0>0$ and $(\pm x_0,\pm y_0)\in C$. We first approximate $\frac{m_{\bm{a}}(z)}{z}$ on $C_1,C_3$.

We use the following approximation lemma.
\begin{lemma}\label{lem:approximate on C1C3}
    Suppose that $F$ has the representation \eqref{eq:Herglotz expression of F}, where $\nu\in\mathcal{E}_p,p\geq\frac{M}{2}$ is a signed measure satisfying $\int_{\mathbb{R}}\frac{1}{1+\lambda^2}\nu(\dd\lambda)<\infty$. Then there is a function $\tilde{F}$, analytic in $D_+$, such that 
    \begin{equation}
        |\tilde{F}(\zeta)-F(-\zeta^2)|\leq \rho_c(\zeta)
    \end{equation}
    on the arcs $C_1,C_3$.
\end{lemma}
\begin{proof}
Fix $c>0$ and, as before, write $H(\zeta)=\cosh{(\zeta^M+K)},\rho_c(\zeta)=H(\zeta)^c$. Choose $b\in D_-$ so that $b^2\notin\mathbb{R}$, and let 
\begin{equation}
    \Phi(\zeta)=\frac{(H(\zeta)H(-\zeta))^{k}}{((\zeta^2-b^2)(\zeta^2-\bar{b}^2))^l}
\end{equation}
where $k,l$ are integers satisfying $2k>c,l\geq2$. The function $\Phi$ has the following properties.\\
(1) Since $\Phi$ is even, $\Phi_\lambda=\Phi(\sqrt{-\lambda})$ is independent of the choice of square root.\\
(2) $\Phi$ is analytic and has neither zeros nor poles in $D_+$.\\
(3) We have $\int_{\mathbb{R}}|\Phi_{\lambda}|\nu(\dd\lambda)<\infty$ whenever $\nu$ satisfies $\int_{\mathbb{R}_-}e^{c|\lambda|^{\frac{M}{2}}}\nu(\dd\lambda)<\infty,\forall c>0$. To see this, split the integral according to the sign of $\lambda$. For $\lambda>0$, the estimate $|\Phi_\lambda|=|\Phi(i\sqrt{\lambda})|\leq C|1+\lambda|^{-2l}$ gives $\int_{\mathbb{R_+}}|\Phi_{\lambda}|\nu(\dd\lambda)<\infty$. For $\lambda<0$, the bound $|\Phi_\lambda|\leq e^{c'|\lambda|^{\frac{M}{2}}}$ gives $\int_{\mathbb{R}_-}|\Phi_{\lambda}|\nu(\dd\lambda)<\infty$.

Write $z=-\zeta^2,t=-\lambda$ and define 
\begin{equation}
    \tilde{F}(\zeta)=\alpha-\beta\zeta^2+\int_{\mathbb{R}}\frac{1-\Phi_\lambda/\Phi(\zeta)}{\lambda+\zeta^2}-\frac{\lambda}{1+\lambda^2}\nu(\dd\lambda).
\end{equation}
To check that the apparent pole at $-\zeta^2=\lambda$ is removable, use the evenness of $\Phi$ to write $\Phi(\zeta)=J(-\zeta^2)$, where $J$ is analytic near the point $-\zeta^2=\lambda$. Then
\begin{equation}
    \frac{1-\Phi_\lambda/\Phi(\zeta)}{\lambda+\zeta^2}=\frac{J(z)-J(\lambda)}{(\lambda-z)J(z)}=\frac{-J'(z)}{J(z)}\in\mathbb{C}.
\end{equation}
At $\lambda=0$, evenness gives $\Phi(\zeta)=\Phi(0)+\Oh(\zeta^2)$, so $1-\Phi_\lambda/\Phi(\zeta)=\Oh(\zeta^2)$. Thus the apparent poles of $\tilde{F}$ are removable.

It remains to verify convergence of the integral. On every compact set $K\subset\subset D_+$, for all sufficiently large $|\lambda|$ we have $\abs{\frac{\Phi_\lambda}{\Phi(\zeta)(\lambda+\zeta^2)}}\leq C_K\frac{|\Phi_\lambda|}{1+|\lambda|}$. The majorant is integrable with respect to $\nu(\dd\lambda)$. Hence $\tilde{F}$ is well defined and analytic in $D_+$.

The error in approximating $F$ by $\tilde{F}$ is
\begin{equation}
    \abs{\tilde{F}(\zeta)-F(-\zeta^2)}\leq\int_{\mathbb{R}}\abs{\frac{\Phi_{\lambda}}{\Phi(\zeta)(\lambda+\zeta^2)}}\nu(\dd\lambda)\leq \frac{C}{|\Im\zeta^2|\Phi(\zeta)}.
\end{equation}

On the arcs $C_1,C_3$, the geometry of the contour gives 
\begin{equation}
    |\Im\zeta^2|=|2ixy|=2|x|\abs{c_0+c_1|x|^{-(M-1)^{-1}}}\geq \frac{C}{2}|x|\geq C|\zeta|
\end{equation} 
which implies 
\begin{equation}
    \abs{\tilde{F}(\zeta)-F(-\zeta^2)}\leq C\la\zeta\ra^{4l-1}|H(\zeta)|^{-2k}.
\end{equation}
On $C_2,C_4$, the bound $\abs{\Im\zeta^2}=|2ixy|\geq 2C|y (y-c_0)^{-(M-1)}|\geq C\la\zeta\ra^{-(M-2)}$ gives 
\begin{equation}
    \abs{\tilde{F}(\zeta)-F(-\zeta^2)}\leq C\la\zeta\ra^{4l+M-2}
\end{equation}
To remove a possible pole of $\frac{m_{\bm{a}}(\zeta)}{\zeta}$ at $\zeta=0$, define 
\begin{equation}
    A_+(\zeta)=\frac{-\hat{m}_+(\zeta)}{\zeta}\left(1-\frac{\Phi(0)}{\Phi(\zeta)}\right),\qquad A_-(\zeta)=\frac{\hat{m}_-(\zeta)}{\zeta}\left(1-\frac{\Phi(0)}{\Phi(\zeta)}\right)
\end{equation}
This function is analytic in $D_+$. The approximation error is
\begin{equation}\begin{aligned}
    \abs{A_+-\frac{m_{\bm{a}}(\zeta)}{\zeta}}&=\abs{\frac{(m_+(-\zeta^2)-\hat{m}_+(\zeta))\Phi(\zeta)+\Phi(0)\hat{m}_+(\zeta)}{\zeta\Phi(\zeta)}}\\
    &\leq C\abs{\la\zeta\ra^{4l-2}H(\zeta)^{-2k}}+|\Phi(0)|\abs{\frac{i\zeta+\Oh(\zeta)}{\zeta}}\abs{\la\zeta\ra^{4l}H(\zeta)^{-2k}}\leq C |\rho_c(\zeta)|^{-1}
\end{aligned}
\end{equation}
For the second term, we use $m_+(\zeta)=\zeta+\Oh(\zeta)$ and $\abs{\tilde{m}_+(\zeta)-m_+(-\zeta^2)}\leq C\la\zeta\ra^{4l-1}|H(\zeta)|^{-2k}$. The same argument gives the corresponding estimate for $A_-(\zeta)-\frac{m_{\bm{a}}(\zeta)}{\zeta}$.
\end{proof}

We now combine the approximations on the four arcs.
\begin{proof}[Proof of Theorem~\ref{Thm: theorem2}, Part I] We first show that $\left(1,\frac{m_{\bm{a}}(\zeta)}{\zeta}\right)\in \mathcal{A}_{N}$. Lemma~\ref{lem:approximate on C1C3} supplies the approximant on $C_1,C_3$. For $C_2$ and $C_4$, we use the expansion of $m_\pm$ in $D_-$. Given $\frac{m_\pm(\zeta)}{\zeta}=1+\sum_{2\leq k\leq N-1}m_{k}\zeta^{-k}+\Oh(\zeta^{-N})$, choose the coefficients $d_{k}$ and set
\begin{equation}
    A_0=1+\sum_{1\leq k\leq N}d_{k}(\zeta-b)^{-k}
\end{equation}
so that $A_0$ matches the stated expansion of $m_\pm$ through the orders $1\leq k\leq N-1$. It follows that
\begin{equation}
    \abs{A_0-\frac{m_\pm(\zeta)}{\zeta}}\leq C|\zeta|^{-N}\leq |\rho_c(\zeta)|^{-1}.
\end{equation}

To combine $A_\pm$ and $A_0$, define 
\begin{equation}
    \eta_j(z)=\frac{\exp{(i^{1-k}\delta z^{4N+1})}}{\sum_{j=1}^{4}\exp{(i^k \delta z^{4 N+1} })}.
\end{equation}
Then $\eta_j(z)\rightarrow 1$ as $z\rightarrow\infty$ along $C_j$, whereas $\eta_j(z)\rightarrow0$ along the other arcs. Set $\tilde{m}(z)=\sum_{j=1}^{4}\eta_j(z)m_j(z)$, where $m_1(\zeta)=A_+(\zeta),m_2(\zeta)=A_0(\zeta),m_3(\zeta)=A_-(\zeta),m_4(\zeta)=A_0(\zeta)$. The zeros of
\begin{equation}
    \sum_{1\leq k\leq4}e^{i^kw}=4\cosh\left(\frac{1+i}{2}w\right)\cosh\left(\frac{1-i}{2}w\right)
\end{equation}
are 
\begin{equation}
    w=\frac{\pi(2k_1+1)}{2}(\pm i\pm1),
\end{equation}
and the corresponding poles of $m_i$ are given by 
\begin{equation}
    \zeta=\left(\frac{\pi(2k_1+1)}{2\delta}\right)^{\frac{1}{4N+1}}\exp\left({i\left({\frac{\pi}{4}+\frac{\pi k_2}{2}+2k_3\pi}\right){(4N+1)}^{-1}}\right)
\end{equation}
Since $\arg(\zeta)\geq\frac{\pi}{4(4N+1)}$, any such point with $|\zeta|>R$ satisfies $\zeta\notin D_+$ once $R$ is sufficiently large. Moreover, $|\zeta|\geq\left(\frac{\pi(2k_1+1)}{2\delta}\right)^{\frac{1}{4N+1}}$. Choosing $\delta$ sufficiently small ensures that $\zeta\notin D_+$ also for $|\zeta|\leq R$. Thus the functions $m_k$ are analytic in $D_+$, and we obtain a bounded analytic function $\tilde{m}(z)$ satisfying 
\begin{equation}
    \sup_{\zeta\in C}|\rho_c(\zeta)|\abs{\frac{m_{\bm{a}}(\zeta)}{\zeta}-\tilde{m}(\zeta)}<\infty,
\end{equation}
which implies $\left(1,\frac{m_{\bm{a}}(\zeta)}{\zeta}\right)\in \mathcal{A}_{N}$.

Set $\bm{b}=\left(1,\frac{2\zeta-m_{\bm{a}}(\zeta)-m_{\bm{a}}(-\zeta)}{m_{\bm{a}}(\zeta)-m_{\bm{a}}(-\zeta)}\right)$. Then $\bm{a}\cdot\bm{b}=\bm{b}\cdot\bm{a}=(1,1)$. If $\Re\zeta>0$, we have $\Re(-\zeta)<0$, and hence 
\begin{equation}
    b_2=\frac{2\zeta+m_{+}(-\zeta^2)-m_{-}(-\zeta^2)}{-m_+(\zeta^2)-m_-(-\zeta^2)}=-2\mathcal{V}_{12}(-\zeta^2)+2\zeta\mathcal{V}_{22}(-\zeta^2)
\end{equation}
We have $\Sigma_{ij}(\dd\lambda)<<\operatorname{tr}\Sigma(\dd\lambda)$, and Lemma~\ref{lem:comparison of sigma and rho} gives $\operatorname{tr}\Sigma\in \mathcal{E}_p$. Consequently, the signed measures representing the entries $\mathcal{V}_{12},\mathcal{V}_{22}$ belong to $\mathcal{E}_p$. Applying Lemma~\ref{lem:approximate on C1C3} gives analytic approximants $\tilde{\mathcal{V}}_{12},\tilde{\mathcal{V}}_{22}$, from which we obtain $\tilde{b}_2$ on $C_1,C_3$. On $C_2,C_4$, we approximate $b_2$ using the expansion of $m_{\bm{a}}$, which gives uniform expansions of $m_{\bm{a},o},m_{\bm{a},e}$ and $b_2$ in $D_-$. We conclude that $\bm{b}\in \mathcal{A}_{N}$ and hence that $\bm{a}\in\mathcal{A}_{N}^{inv}$.

By uniqueness of analytic continuation, the $m$-function extends from $D_{-}$ to $\mathbb{C}\setminus(\mathbb{R}\cup i\mathbb{R})$. The argument of \cite[Proposition~4.5]{Kotani2024} then gives $\bm{a}\in \mathcal{A}_{N,+}^{inv}$.

The preceding construction uses the following three conditions on $m$.\\
$(M.1)$ The function $m(z)$ satisfies the reality and Herglotz conditions
\begin{equation}
    m(\bar{z})=\overline{m(z)},\qquad \frac{\Im{m(z)}}{\Im{z}}>0\text{ on }\mathbb{C}\setminus(\mathbb{R}\cup i\mathbb{R}).
\end{equation}
$(M.2)$ The function $m(z)$ has the asymptotic expansion
\begin{equation}
    m(z)=z+\sum_{1\leq k\leq N}m_k z^{-k}+O(z^{-N+1}) \text{ on }D_-.
\end{equation}
$(M.3)$ The spectral measures $\sigma_\pm$ satisfy the moment condition
\begin{equation}
    \int_{\mathbb{R}_-}e^{c|\lambda|^{\frac{M}{2}}}\sigma_{\pm}(d\lambda)<\infty \text{ for any }c>0.
\end{equation}
Under these assumptions, Lemmas~\ref{lem:comparison of sigma and rho}, \ref{lem:approximate on C1C3}, and~\ref{lem:herglotz for mga} give the first part of Theorem~\ref{Thm: theorem2}.
    
\end{proof}

Thus the symbol associated with the initial data belongs to $\mathcal{A}_{N,+}^{inv}$. 

\subsection{Preservation of the spectral conditions and the KdV flow}
We prove that the action $m_{\bm{a}}\mapsto m_{g\bm{a}}$ preserves conditions $(M.1),(M.2),(M.3)$ and reconstruct the solution $q(x,t)$ of the KdV equation.
\begin{lemma}\label{lem:herglotz for mga}
    Let $\bm{a}\in\mathcal{A}_{N,+}^{inv}$ have an $m$-function satisfying $(M.1)$, and let $g=re^h\in \Gamma_{m}^{real}$ with $\deg r=0$. If $g^{-1}T(g\bm{a})$ is invertible on $H_{N,c}(D_+)$, then $m_{g\bm{a}}$ is well defined on $D_-$ and extends to $\mathbb{C}\setminus(\mathbb{R}\cup i\mathbb{R})$, where it satisfies $(M.1)$. The parameters $M,N,c$ are subject to the conditions
    \begin{equation}
        M\ge m,\qquad N\ge m+1,\qquad\begin{cases}
        c>0,&M>m,\\
        c>c_2,&M=m,
        \end{cases}
    \end{equation}
    where $c_2$ is given by 
    \begin{equation}
        |r_n(z)|+|r_n(z)^{-1}|+|g(z)|+|g(z)^{-1}| \leq C e^{c_2|\Re z|^m}, \qquad z\in C_x.
    \end{equation}
\end{lemma}
\begin{proof}
    It suffices to treat $g=e^h$ with $\operatorname{deg}h\leq m$. Set 
    \begin{equation}
        B_n=r_n^{-1}T(r_n\bm{a}),\qquad B=g^{-1}T(g\bm{a}).
    \end{equation}
    where the approximants from Lemma~\ref{rneh} satisfy $r_n\rightarrow g$ in the sense of \eqref{in the sense of}. Applying \eqref{eq:estimate of rn x}, \eqref{eq: estimate of rn y}, and the dominated convergence theorem to $d_c(g,r_n)^2$ gives 
    \begin{equation}
        \norm{B-B_n}_{H.S.}\rightarrow0
    \end{equation}
    as $n\rightarrow\infty$. Consequently,
    \begin{equation}
        \norm{B^{-1}-B_n^{-1}}_{H.S.}\leq\norm{B^{-1}}_{\mathcal{B}(H)}\norm{B_n-B}_{H.S.}\norm{B_n^{-1}}_{\mathcal{B}(H)}\leq 2\norm{B^{-1}}_{\mathcal{B}(H)}\norm{B_n-B}_{H.S.}\rightarrow0.
    \end{equation}
    This estimate requires $N\geq m+1$.
    
    Define 
    \begin{equation}\left\{\begin{aligned}
        B_n^{-1}(r_n^{-1}1)=u_n,\qquad B^{-1}(g^{-1}1)=u,\\
        B_n^{-1}(r_n^{-1}z)=v_n,\qquad B^{-1}(g^{-1}z)=v.
    \end{aligned}\right.
    \end{equation}
    and 
    \begin{equation}\left\{\begin{aligned}
        \varphi_{n}=r_n\bm{a}u_n-1, \qquad\varphi=g\bm{a}u-1,\\
        \psi_{n}=r_n\bm{a}v_n-z, \qquad\psi=g\bm{a}v-z.
    \end{aligned}\right.
    \end{equation}
    Another application of dominated convergence gives $r_n\rightarrow g$ and $r_n^{-1}\rightarrow g^{-1}$ in $H_{N,c}(D_+)$. Hence 
    \begin{equation}
        \norm{u_n-u}_{H_{N,c}(D_+)}\leq\norm{B_n^{-1}(r_n^{-1}1)-B^{-1}(r_n^{-1}1)}_{H_{N,c}(D_+)}+\norm{B^{-1}(r_n^{-1}1)-B^{-1}(g^{-1}1)}_{H_{N,c}(D_+)}\rightarrow0
    \end{equation}
    and likewise $v_n\rightarrow v$ in $H_{N,c}(D_+)$. Choose an approximant $\tilde{\bm{a}}$ and fix $c>c_2$. Then 
    \begin{equation}\begin{aligned}
        \norm{r_n(\bm{a}-\tilde{\bm{a}})u_n-g(\bm{a}-\tilde{\bm{a}})u}_{L^2(C)}&\leq\norm{r_n(\bm{a}-\tilde{\bm{a}})(u_n-u)}_{L^2(C)}+\norm{(r_n-g)(\bm{a}-\tilde{\bm{a}})u}_{L^2(C)}\\
    \end{aligned}
    \end{equation}
    The first term is bounded by $C\|u_n-u\|_{H_{N,c}(D_+)}\rightarrow0$, and the second tends to zero pointwise. Dominated convergence therefore gives 
    \begin{equation}
        \norm{r_n(\bm{a}-\tilde{\bm{a}})u_n-g(\bm{a}-\tilde{\bm{a}})u}_{L^2(C)}\rightarrow0.
    \end{equation}
    Since $\mathfrak{p}_-:L^2(C)\rightarrow H(D_-)$ is a bounded operator and 
    \begin{equation}
        \mathfrak{p}_{-}(r_n\bm{a}u_n)=\mathfrak{p}_{-}(r_n(\bm{a}-\tilde{\bm{a}})u_n)=\varphi_n,\qquad \mathfrak{p}_{-}(g\bm{a}u)=\mathfrak{p}_{-}(g(\bm{a}-\tilde{\bm{a}})u)=\varphi,
    \end{equation}
    we obtain 
    \begin{equation}
        \norm{\varphi_n-\varphi}_{H(D_-)}\rightarrow0
    \end{equation}
    as $n\rightarrow\infty$. 
    
    On each compact set $K\subset D_-$, the Cauchy--Schwarz inequality gives
    \begin{equation}
        \sup_{z\in K}|\varphi_n(z)-\varphi(z)|\leq\frac1{2\pi}\sup_{z\in K}\left(\int_C\frac{|\dd\lambda|}{|\lambda-z|^2}\right)^{1/2}\norm{r_n(\bm{a}-\tilde{\bm{a}})u_n-g(\bm{a}-\tilde{\bm{a}})u}_{L^2(C)}\rightarrow0.
    \end{equation}
    Hence $\varphi_n\rightarrow\varphi$ locally uniformly in $D_-$. The same argument gives $\psi_n\rightarrow\psi$ and, on taking the quotient, $m_{r_n\bm{a}}\rightarrow m_{g\bm{a}}$. 

    We define
    \begin{equation}
        m_{+,r_n\bm{a}}(z)=-m_{r_n\bm{a}}(\sqrt{-z}),\qquad m_{-,r_n\bm{a}}(z)=m_{r_{n}\bm{a}}(-\sqrt-z).
    \end{equation} 
    By assumption, $m_{+,\bm{a}},m_{-,\bm{a}}$ are Herglotz functions. Since $r_n=r_{\zeta_{m_1}}r^{-1}_{\eta_1}\cdots r^{-1}_{\eta_{m_2}}$ with $\zeta_j,\eta_j\in D_-$, the Darboux transformation in \cite[Lemma 7.2]{Kotani2024} gives Herglotz extensions of $m_{+,r_n\bm{a}},m_{-,r_n\bm{a}}$. Their locally uniform limits $m_{+,g\bm{a}},m_{-,g\bm{a}}$ are likewise Herglotz functions. Thus $m_{g\bm{a}}$ satisfies $(M.1)$ for every $g=re^h\in \Gamma_{m}^{real}$ with $\deg r=0$.
\end{proof}

With the characteristic functions $\varphi_{g\bm{a}},\psi_{g\bm{a}}$ and the associated function $m_{g\bm{a}}$ in hand, define 
\begin{equation}
    (K(g)q)(x)=\mathcal{R}(g\bm{a})=-2\partial_x\kappa_1(e_xg\bm{a}).
\end{equation}

The expansion required by $(M.2)$ for $m_{g\bm{a}}$ follows from \cite[Proposition 4.7]{Kotani2024}.
\begin{lemma}\label{lem: curve C and sigma C}
    If ${\bm{a}}\in \mathcal{A}_{N,+}^{inv}(C)$, then $m_{g\bm{a}}$ satisfies $(M.2)$ to order $N-\frac{m+1}{2}$ on $\sigma D_-$ for $\sigma>1$.
\end{lemma}

Lemma~\ref{lem:herglotz for mga} gives representing measures for the Herglotz functions $m_{+,g\bm{a}},m_{-,g\bm{a}}$. We verify $(M.3)$ for these measures using the transfer matrix. At this stage the limit-point property for the potential $K(g)q$ is still to be proved, so uniqueness of the spectral measures is not yet available. 

Recall that $S=m_++m_-$ and 
\begin{equation}
    \mathcal{V}=\frac{1}{S}\left(\begin{array}{cc}
         m_+m_-&(m_+-m_-)/2  \\
         (m_+-m_-)/2&-1 
    \end{array}\right).
\end{equation}

For the transfer matrix, we use the notation of \cite[Chapter~6]{Kotani2024}. Set
\begin{equation} \left\{\begin{aligned}
    \Pi_{\bm{a}}(\zeta)&=\left(\begin{array}{cc}
         \psi_{\bm{a},e}(\zeta)+\kappa_{1}(\bm{a})(1+\varphi_{\bm{a}}(\zeta))&\zeta+\psi_{\bm{a},o}(\zeta)+\kappa_1(\bm{a})\varphi_{\bm{a},o}(\zeta)  \\
         1+\varphi_{\bm{a},e}(\zeta)& \varphi_{\bm{a},o}(\zeta)
    \end{array}\right)\\
    G(\zeta)&=\left(\begin{array}{cc}
         g_e(\zeta)&g_o(\zeta)  \\
         g_o(\zeta)&g_e(\zeta) 
    \end{array}\right)\text{ for }g=re^h\in\Gamma_{m}\text{ with }\deg r=0
\end{aligned}\right.
\end{equation}
Since $\det(G(\zeta))\neq0,\det(\Pi_{\bm{a}(\zeta)})\neq0$ on $D_-$, define
\begin{equation}
    T_{\bm{a}}(-\zeta^2,g)=\Pi_{g\bm{a}}(\zeta)G(\zeta)^{-1}\Pi_{\bm{a}}(\zeta)^{-1},
\end{equation}
so that
\begin{equation}
    \left\{\begin{aligned}
        &T_{\bm{a}}(-\zeta^2,g)\bm{\pi}_a(\zeta)=g(\zeta)^{-1}\bm{\pi}_{g\bm{a}}(\zeta)\\
        &\bm{\pi}_{\bm{a}}(\zeta)=\left(\begin{matrix}
            \zeta+\psi_{\bm{a}}(\zeta)+\kappa_1(\bm{a})(1+\varphi_{\bm{a}}(\zeta))\\
            1+\varphi_{\bm{a}}(\zeta)
        \end{matrix}\right)
    \end{aligned}\right.,
\end{equation}
where $z=-\zeta^2$ is the spectral parameter.

Since $m_{a}(\zeta)=\frac{(\bm{\pi}_a)_1}{(\bm{\pi}_a)_2}$, the preceding identity gives 
\begin{equation}\label{eq:transfer matrix of ma}
    m_{g\bm{a}}(\zeta)=T_{\bm{a}}(-\zeta^2,g)\cdot m_{\bm{a}}(\zeta)
\end{equation}
where the fractional linear action is given by $T\cdot w=\frac{T_{11}w+T_{12}}{T_{21}m+T_{22}}$. The sign conventions for the two branches are
\begin{equation}
    m_+^q(z)=-m_{\bm{a}}(\sqrt{-z}),\qquad m_-^q(z)=m_{\bm{a}}(-\sqrt{-z}).
\end{equation}
Thus $u,v$ transform under the same matrix according to
\begin{equation}
    u_g=T\cdot u,\qquad v_g=T\cdot v,\qquad T=T_q(z,g).
\end{equation}

Define
\begin{equation}
    g_s=e^{sh},\qquad q_s=K(g_s)q,\qquad0\leq s\leq1,\qquad T_s(z)=T_q(z,g_s)
\end{equation}
and assume the following.\\
$(T.1)$ The functions $m_\pm^{q_s}$ are Herglotz functions with spectral measures $\sigma_\pm^{q_s}$.\\
$(T.2)$ The transfer matrix implements the action $K(g)$ through \eqref{eq:transfer matrix of ma}. For $g_1,g_2\in\Gamma_m$ with the rational function part of degree 0, it satisfies the cocycle identity
\begin{equation}
    T_{\bm{a}}(z,g_1g_2)=T_{g_1\bm{a}}(z,g_2)T_{\bm{a}}(z,g_1).
\end{equation}
Moreover, for each $s$, the function $m_s$ has the expansion
\begin{equation}
    m_s(\zeta)=\zeta+\sum_{k=1}^{N}m_{k}(s)\zeta^{-k}+\Oh(\zeta^{-N-1})
\end{equation}
For fixed $g=e^h$, the matrix-valued function $T_{\bm{a}}(\cdot,e^h)$ is entire and satisfies
\begin{equation}
    T_{\bm{a}}(\bar{z},e^h)=\overline{T_{\bm{a}}(z,e^h)},\qquad \det T_{\bm{a}}(z,e^h)=1
\end{equation}

Matrix multiplication and the identity $\det T=1$ give
\begin{equation}
    e(T\cdot w)=\frac{Te(w)}{T_{21}w+T_{22}},T\cdot u-T\cdot v=\frac{(T_{11}T_{22}-T_{12}T_{21})(u-v)}{(T_{21}u+T_{22})(T_{21}v+T_{22})}=\frac{(u-v)}{(T_{21}u+T_{22})(T_{21}v+T_{22})}
\end{equation}
Substitution gives the transformation law for $\mathcal{V}_{K(g)q}$,
\begin{equation}\label{eq:V K(g)q,V q}\begin{aligned}
    \mathcal{V}_{K(g)q}&=\frac{e(T\cdot u)e(T\cdot v)^T+e(T\cdot v)e(T\cdot u)^T}{2(T\cdot u-T\cdot v)}\\
    &=\frac{T(e(u)e(v)+e(v)e(u))T^T}{(T_{21}u+T_{22})(T_{21}v+T_{22})}\frac{(T_{21}u+T_{22})(T_{21}v+T_{22})}{2(u-v)}\\
    &=T\mathcal{V}_qT^T
\end{aligned}
\end{equation}
\begin{lemma}
    The corresponding matrix spectral measures satisfy
    \begin{equation}\label{eq:equivalence of Sigma}
        \Sigma_{K(g)q}(\dd\lambda)=T(\lambda)\Sigma_q(\dd\lambda)T(\lambda)^T
    \end{equation}
\end{lemma}
\begin{proof}
    We apply Stieltjes inversion against a test function. Fix $R>0$ and choose $\varphi\in C_c(\mathbb{R})$ with $\operatorname{supp}(\varphi)\subset K_R=[-R,R]$. The Herglotz representation gives
    \begin{equation}
        \mathcal{V}_q(z)=\mathcal{V}_R(z)+\int_{K_R}\frac{\Sigma_q(\dd \lambda)}{\lambda-z}
    \end{equation}
    with 
    \begin{equation}
        \mathcal{V}_R=A+Bz+\int_{K_R}\frac{\lambda}{1+\lambda^2}\Sigma_q(\dd\lambda)+\int_{\mathbb{R}\setminus K_R}\left(\frac{1}{\lambda-z}-\frac{\lambda}{1+\lambda^2}\right)\Sigma_q(\dd\lambda)
    \end{equation}
    The part $\mathcal{V}_R$ extends analytically across the interior of $K_R$ and has the symmetry
    \begin{equation}
        \mathcal{V}_R(x)^*=\mathcal{V}_R(x),\qquad \forall x\in K_R
    \end{equation}
    
    By \eqref{eq:V K(g)q,V q},
    \begin{equation}\label{eq: def of F R}
        \mathcal{V}_{K(g)q}(z)=\int_{K_R}\frac{T(\lambda)\Sigma_q(\dd\lambda)T(\lambda)^T}{\lambda-z}+\mathcal{F}_R(z)
    \end{equation}
    where the remainder is
    \begin{equation}
           \mathcal{F}_R(z)=T(z)\mathcal{V}_R(z)T(z)^T+\int_{K_R}\frac{T(z)\Sigma_q(\dd\lambda)T(z)^T-T(\lambda)\Sigma_q(\dd\lambda)T(\lambda)^T}{\lambda-z}.
    \end{equation}
    The integral in the remainder splits as 
    \begin{equation}\begin{aligned}
        \int_{K_R}T(z)\Sigma_q(\dd\lambda)\frac{(T(z)-T(\lambda))^T}{\lambda-z}+\frac{T(z)-T(\lambda)}{\lambda-z}\Sigma_q(\dd\lambda)T(\lambda)^T.
    \end{aligned}
    \end{equation}
    For fixed $\lambda$, the divided difference $\frac{T(z)-T(\lambda)}{\lambda-z}$ is analytic in $z$ after removal of the singularity at $z=\lambda$. Thus $\mathcal{F}_R$ extends analytically across the interior of $K_R$ and satisfies
    \begin{equation}
        \mathcal{F}_R(x)^*=\mathcal{F}_R(x),\qquad \forall x\in K_R
    \end{equation}
    In particular, on compact subintervals of the interior we have 
    \begin{equation}
        \Im(\mathcal{F}_R(x+iy))=\Oh_R (y)
    \end{equation}
    Taking $R$ large enough that the support of the test function lies in the interior and applying Stieltjes inversion to \eqref{eq: def of F R} gives
    \begin{equation}
        \int_{\mathbb{R}}\varphi(\lambda)\Sigma_{K(g)q}(\dd\lambda)=\lim_{y\rightarrow0+}\frac{1}{\pi}\int_{\mathbb{R}}\varphi(x)\Im\mathcal{V}_{K(g)q}(x+iy)\dd x=\int_{K_R}\varphi(\lambda)T(\lambda)\Sigma_q(\dd\lambda)T(\lambda)^T.
    \end{equation}
    Since $R>0,\varphi\in C_c(K_R)$ were arbitrary, we conclude that 
    \begin{equation}
        \Sigma_{K(g)q}(\dd\lambda)=T(\lambda)\Sigma_q(\dd\lambda)T(\lambda)^T
    \end{equation}

\end{proof}
Since $\det T=1$, we have $\|T^{-1}\|=\|T\|$. Write $\Sigma_q=H_q\rho_q$, where the matrix density satisfies $H_q\geq0, \operatorname{tr}H_q=1$ almost everywhere with respect to $\rho_q$. For $\lambda\in \mathbb{R}$, the matrix $T$ is real, so
\begin{equation}
    \|T\|^{-2}\leq\operatorname{tr}(TH_qT^T)\leq\|T\|^2.
\end{equation}
Hence \eqref{eq:equivalence of Sigma} gives 
\begin{equation}\label{eq:estimate of rho K(g)q}
    \|T(\lambda)\|^{-2}\rho_q(\dd\lambda)\leq \rho_{K(g)q}(\dd\lambda)\leq\|T(\lambda)\|^{2}\rho_q(\dd\lambda).
\end{equation}
\begin{lemma}
    Under assumption $(T.2)$, the transfer matrix satisfies 
    \begin{equation}\label{eq:estimate of Ts}
        \|T_q(-R,e^{sh})\|^2\leq C(1+R)e^{CR^{d/2}},\qquad R\geq1,\ 0\leq s\leq1
    \end{equation}
\end{lemma}

\begin{proof}
    Write $T_s=T_q(z,e^{sh})$. Equations (6.5), (6.11), and (6.12) of \cite{Kotani2024} give
    \begin{equation}
        T_{s+\tau}=T_{q_s}(z,e^{\tau h})T_s,\qquad A_{q_s}(z,h)=\partial_{\tau}T_{q_s}(z,e^{\tau h})|_{\tau=0}
    \end{equation}
    Differentiating with respect to $\tau$ at zero gives
    \begin{equation}\label{eq:ODE of T}
        \partial_s T_s=A_{q_s}(z,h) T_s,\ T_0=I
    \end{equation}
    Write $m_s=m_{e,s}+m_{o,s}$ for the decomposition into even and odd parts. \cite[Lemma~6.2]{Kotani2024} gives
    \begin{equation}
        A_{q_s}(-\zeta^2,h)=\mathfrak{p}_+\left[\frac{h}{m_{o,s}}\left(\begin{array}{cc}
             -m_{e,s}&m_{e,s}^2-m_{o,s}^2  \\
             -1& m_{e,s}
        \end{array}\right)\right]+\left(\begin{array}{cc}
             0& c_s(h) \\
             0&0 
        \end{array}\right)
    \end{equation}
    with
    \begin{equation}
        c_{s}(h)=\lim_{\zeta\rightarrow\infty}\left[h+\mathfrak{p}_+\left(\frac{m_{e,s}h}{m_{o,s}}\right)-m_{s}\mathfrak{p}_+\left(\frac{h}{m_{o,s}}\right)\right]
    \end{equation}

    The expansion of $m_s$ gives $m_{o,s}=\zeta+\Oh(\zeta^{-1}),m_{e,s}=\Oh(\zeta^{-2})$. Writing $\deg h(\zeta)=d=2k+1$, we obtain the following bounds on the degrees of the entries of $A$,
    \begin{equation}
        \deg A_{12}\leq k+1,\qquad\deg A_{21}\leq k,\qquad\deg A_{11},\deg A_{22}\leq k-1
    \end{equation}
    regarded as polynomials in $z=-\zeta^2$. Their coefficients are uniformly bounded by continuity in $s\in[0,1]$.

    Write $R=\kappa^2,\kappa\geq1$ and set
    \begin{equation}
        D_{\kappa}=\operatorname{diag}(\kappa^{-1},1),\qquad \tilde{T}_s=D\kappa T_sD_\kappa^{-1},\qquad \tilde{A}_s=D\kappa A_sD_\kappa^{-1}.
    \end{equation}
    The degree bounds for the entries of $A$ give 
    \begin{equation}
        |A_{12}|\leq C\kappa^{d+1},\qquad|A_{12}|\leq C\kappa^{d-1},\qquad|A_{12}|\leq C\kappa^{d-3}
    \end{equation}
    Conjugation to form $D_\kappa A D_\kappa^{-1}$ affects only the off-diagonal entries $A_{12}$ and $A_{21}$. Thus
    \begin{equation}
        \|\tilde{A}_s(-\kappa^2,h)\|\leq C\kappa^d.
    \end{equation}
    After conjugation, \eqref{eq:ODE of T} takes the integral form $\tilde{T}_t=I+\int_{0}^{t}\tilde{A}_{q_s}\tilde{T}_s\dd s$. Consequently,
    \begin{equation}
        \|\tilde{T}_t\|\leq I+\int_{0}^{t}\|\tilde{A}_{q_s}\|\|\tilde{T}_s\|\dd s=:w(t).
    \end{equation}
    Differentiation gives
    \begin{equation}
        w'(t)=\|\tilde{A}_{q_t}\|\|\tilde{T}_t\|\leq\|\tilde{A}_{q_t}\| w(t),
    \end{equation}
    Gr\"onwall's inequality now gives
    \begin{equation}
        \|\tilde{T}_t\|\leq w(t)\leq w(0)\exp\left(\int_{0}^{t}\|\tilde{A}_{q_s}\|\dd s\right)\leq e^{C\kappa^d}.
    \end{equation}
    Since $\|D_\kappa\|\|D_\kappa^{-1}\|=\kappa$, we have 
    \begin{equation}
        \|T_q(-R,e^{sh})\|^2\leq \left(\kappa e^{C\kappa^d}\right)^2\leq C(1+R)e^{C'R^{d/2}}.
    \end{equation}
\end{proof}

We can now prove preservation of the negative spectral moments.
\begin{lemma}\label{lem:(M.3)}
    Under assumptions $(T.1),(T.2)$, the KdV flow preserves $(M.3)$ in the sense that
    \begin{equation}     \sigma_{+}^q+\sigma_{-}^q\in\mathcal{E}_p\iff\sigma_{+}^{K(g)q}+\sigma_{-}^{K(g)q}\in\mathcal{E}_p
    \end{equation}
\end{lemma}
\begin{proof}
    By Lemma~\ref{lem:comparison of sigma and rho}, $\sigma_{+}^q+\sigma_{-}^{q}\in\mathcal{E}_p$ implies $\rho_{q}\in\mathcal{E}_p$. Combining \eqref{eq:estimate of rho K(g)q} and \eqref{eq:estimate of Ts} gives
    \begin{equation}\begin{aligned}
        \int_{-\infty}^{-1}e^{c|\lambda|^p}\rho_{K(g)q}(\dd\lambda)&\leq C\int_{-\infty}^{-1}(1+|\lambda|)e^{c|\lambda|^p+C|\lambda|^{d/2}}\rho_q(\dd\lambda)\\
        &\leq C'\int_{-\infty}^{-1}e^{(c+C+1)|\lambda|^p}<\infty
    \end{aligned}
    \end{equation}
    for $p\geq d/2$, where we have used $1+R\leq C_p e^{R^p}$. Since $\rho([-1,0])<\infty$, it follows that $\rho_{K(g)q}\in\mathcal{E}_p$. Another application of Lemma~\ref{lem:comparison of sigma and rho} gives $\sigma_{+}^{K(g)q}+\sigma_{-}^{K(g)q}\in\mathcal{E}_p$. For the converse, use the other inequality in \eqref{eq:estimate of rho K(g)q},
    \begin{equation}
        \rho_q(\dd\lambda)\leq \|T(\lambda)\|^{2}\rho_{K(g)q}(\dd\lambda)
    \end{equation}
    and apply the same argument.
\end{proof}
Lemma~\ref{lem:herglotz for mga} gives $(T.1)$, and the transfer identities give $(T.2)$ for $m_{\pm}^{q_s}$. Hence the KdV evolution constructed above preserves $(M.3)$. 

We next deduce the limit-point property from $(M.3)$ and reconstruct the KdV solution from a function $m(z)$ satisfying $(M.1),(M.2),(M.3)$.

\begin{lemma}\label{lem:LP}
    Let $q\in L^{1}_{loc}(\mathbb{R})$ be real. If a half-line spectral measure $\sigma_{\pm}$ satisfies $(M.3)$, then the corresponding endpoint $\pm\infty$ is limit point.
\end{lemma}
\begin{proof}
    We give the proof for $+\infty$. Suppose to the contrary that $+\infty$ is limit circle. Then, for every $z\in\mathbb{C}$, we have $\dim(W_z)=2$, where $W_z$ denotes the space of square-integrable solutions of
    \begin{equation}
        -f''(x,z)+q(x)f(x,z)=zf(x,z),\ x\geq0
    \end{equation}
    Choose zero-energy solutions $u,v$ with 
    \begin{equation}
        \left\{\begin{array}{ccc}
            -u''(x)+q(x)u(x)=0, &u(0)=0, &u'(0)=1,\\
            -v''(x)+q(x)v(x)=0, &v(0)=1, &v'(0)=0,
        \end{array}\right.
    \end{equation}
    Since $+\infty$ is limit circle, $u,v\in L^2([0,\infty))$. Set 
    \begin{equation}
        h(x)=\sqrt{v^2(x)+u^2(x)}>|u(x)|,|v(x)|,\qquad\ A=\int_{0}^{\infty}h^2(x)\dd x<\infty.
    \end{equation}
    At the regular endpoint $0$, take $s(x,z)$ to be the solution of
    \begin{equation}
        -s''(x,z)+q(x)s(x,z)=zs(x,z),\qquad s(0,z)=0,\qquad s'(0,z)=1.
    \end{equation}
    Since $s(x,0)=u(x)$, variation of parameters gives
    \begin{equation}
        s(x,z)=u(x)+z\int_{0}^{x}\left(u(y)v(x)-u(x)v(y)\right)s(y,z)\dd y.
    \end{equation}
    The bound $\abs{u(y)v(x)-u(x)v(y)}\leq 2h(x)h(y)$ gives
    \begin{equation}
        \abs{\frac{s(x,z)}{h(x)}}\leq1+2|z|\int_{0}^{x}h^2(y)\abs{\frac{s(y,z)}{h(y)}}\dd y.
    \end{equation}
    We use Gr\"onwall's inequality in the form
    \begin{equation}
        w(x)\leq\alpha+\int_{a}^{x}\beta(y)w(y)\dd y\implies w(x)\leq\alpha\exp\left(\int_{a}^{x}\beta(y)\dd y\right)
    \end{equation}
    Here we take 
    \begin{equation}
        w(x)=\abs{\frac{s(x,z)}{h(x)}},\qquad \beta(y)=2|z|h^2(y),\ a=0,\ \alpha=1
    \end{equation}
    to obtain 
    \begin{equation}
        |s(x,z)|\leq h(x)\exp\left(2|z|\int_{0}^{x}h^2(y)\dd y\right),
    \end{equation}
    and hence
    \begin{equation}
        \|s(x,z)\|^2_{L^2(\mathbb{R})}=\int_{0}^{\infty}\abs{s(x,z)}^2\dd x\leq Ae^{2A|z|}
    \end{equation}

    Under the limit-circle hypothesis at $+\infty$, the half-line realization of $L_q$ has pure point spectrum with a sequence of eigenvalues tending to $-\infty$ (see \cite{GesztesySimon2000} and \cite{Matsuda1966}). Denote the corresponding eigenpairs by $(\lambda_n,s_n) $, with $\lambda_n\rightarrow-\infty.$

    In the normalization above, the mass at each eigenvalue is given by
    \begin{equation}
        \sigma_{+}(\{\lambda_n\})=\frac{1}{\|s_n\|_{L^2([0,\infty))}^2}
    \end{equation}
    by Lemma~9.13 and the example following it in \cite{Teschl2014}. Hence
    \begin{equation}
        \int_{-\infty}^0e^{c|\lambda|^\frac{M}{2}}\sigma_{+}(\dd\lambda)\geq\sum_{n=1}^{\infty}e^{c|\lambda_n|^\frac{M}{2}}A^{-1}e^{-2A|\lambda_n|}.
    \end{equation}
    The right-hand side equals $\infty$ if $M>2$, or if $M=2,c-2A>0$. This contradicts $(M.3)$.
\end{proof}
With the choice of $c$ above, the same argument gives the implication $\int_{\mathbb{R}_-}e^{c|\lambda|}\sigma(\dd\lambda)<\infty\implies LP$. 

We now construct the KdV flow.

\begin{proof}[Proof of Theorem~\ref{Thm: theorem2}, Part II.] Let $q_0\in \mathcal{Q}_{N}^{(m)}$. There is a symbol $\bm{a}_0\in\mathcal{A}_{N,+}^{inv}$ such that 
    \begin{equation}
        q_0=\mathcal{R}(\bm{a}_0)
    \end{equation}
    By Lemma~\ref{lem:tau ga(g)>0}, $\tau^{(2)}_{\bm{a}_0}(e_xg)>0$ for all $g=re^h\in\Gamma_{m}^{real},\deg r=0,x\in\mathbb{R}$. Thus $(e_xg)^{-1}T(e_xg\bm{a}_0)$ is invertible, and we can define the two characteristic functions and the coefficient
    \begin{equation}
        \varphi_{e_xg\bm{a}_0},\qquad\psi_{e_xg\bm{a}_0},\qquad\kappa_1(e_xg\bm{a}_0).
    \end{equation}
     The reconstructed potential is given by 
    \begin{equation}
        (K(g)q_0)(x)=\mathcal{R}(g\bm{a}_0)=-2\partial_x\kappa_1(e_xg\bm{a}_0).
    \end{equation}
    By Lemmas~\ref{lem:herglotz for mga}, \ref{lem: curve C and sigma C}, and~\ref{lem:(M.3)}, the action $K(g)$ preserves $(M.1),(M.2),(M.3)$, with the order in $(M.2)$ reduced to $N-\frac{m+1}{2}$.

\end{proof}

\section{Sufficient conditions on the potential}
We give conditions on the potential that imply the required Weyl asymptotics and spectral moment bounds. For the asymptotics we use a contraction argument for a Volterra equation.

\subsection{Asymptotic expansion of the Weyl function}

We establish a uniform expansion of the Weyl $m$-function for potentials that may be unbounded below. Recall the region in the $z$-plane
\begin{equation}
    \mathcal{D}=\{z:\Re z\geq k_0\geq C_0+1, |\Im z|\geq C(\Re z)^{-M/2}\}.
\end{equation}
Here $M\geq0$, and $N\geq3$ denotes the order of the expansion. We impose the following conditions on $q$.
\begin{lemma}\label{lem: expansion}
    Let $r\geq M+N+2$ be an even integer, and suppose that $q\in C^r(\mathbb{R})$ satisfies the following conditions.
\begin{enumerate}[{$(A)$}]
    \item There is a constant $C_0>0$ such that $q(x)\leq C_0$. The half-line operators $L_q^{\pm}:=(-\partial_x^2+q)|_{\mathbb{R}_{\pm}}$ are limit point at their respective endpoints $\pm\infty$.
    \item With $Q(x)=1+C_0-q(x)\geq1$ and $a=\frac{1}{2}(1-\frac{M+N+2}{r})\in[0,\frac{1}{2})$, we assume the bounds for any $1\leq j\leq r$
    \begin{equation}
        \abs{q^{(j)}(x)}\leq C_rQ(x)^{1+ja}.
    \end{equation}
\end{enumerate}

\end{lemma}

\subsubsection{Derivative estimates}

We prove the expansion of $m_+$ in $\mathcal{D}\cap\{\Re z>\Im z>0\}$. Since $\Re z\geq C_0+1\geq q+1$, we have $|z-q(x)|\geq1$. Define 
\begin{equation}
    p(x,z)=\sqrt{z-q(x)},\qquad \rho(x)=Q(x)^a,\qquad h(x,z)=\frac{\rho(x)}{|p(x,z)|},\qquad D=p(x,z)^{-1}\partial_x
\end{equation}
We have $|p|^2\asymp \Re z+Q$ in the region $\{(x,z):x\in\mathbb{R},\Re z>\Im z\geq C_0+1\}$.

If $Q(x)\leq \Re z$, then
\begin{equation}
    h=\frac{\rho}{|\sqrt{z-q}|}\leq C\frac{Q^{a}}{|z|^{\frac{1}{2}}}\leq C'|z|^{a-\frac{1}{2}},
\end{equation}
whereas if $Q(x)\geq \Re z$, then
\begin{equation}
    h=\frac{\rho}{|\sqrt{z-q}|}\leq C\frac{Q(x)^a}{Q(x)^{\frac{1}{2}}}\leq C' |z|^{a-\frac{1}{2}}.
\end{equation}
Thus $\sup_{x\geq0}h(x,z)\rightarrow0$ as $|\Re z|\rightarrow\infty$. Choosing $k_0$ sufficiently large gives $\sup_{x\geq0} h(x,z)\leq\frac{1}{2}$ throughout $\mathcal{D}$.

We first estimate the derivatives of $Q,\rho$ and $h$.

Write $f\in\Sh_m^n$ if $f\in C^n$ and $|D^j(f)|\leq C_jh^{j+m}$ for $0\leq j\leq n$. Likewise, write $f\in\Th_m^n$ if $f\in C^n$ and $|D^j(f)|\leq C_jh^{j+m}|p|^{-2}$ for $0\leq j\leq n$. The Leibniz rule for $D$ gives
\begin{equation}\label{eq:relationship between S,T}
        \Sh_{m_1}^n\Sh_{m_2}^n\subset \Sh_{m_1+m_2}^n,\Sh_{m_1}^n\Th_{m_2}^n\subset \Th_{m_1+m_2}^n, D\Sh_{m}^n\subset \Sh_{m+1}^{n-1},D\Th_{m}^n\subset \Th_{m+1}^{n-1}.
\end{equation}

\subsubsection{The Liouville transformation and higher-order approximate frequencies}

The spectral equation is
\begin{equation}\label{eq:schrodinger}
  y''+p(x,z)^2y=0.
\end{equation}
Make the Liouville substitution
\begin{equation}\label{eq:liouville}
  y=p^{-1/2}u.
\end{equation}
Direct differentiation gives
\begin{equation}\label{eq:liouville-equation}
  D^2u+(1+V)u=0,
\end{equation}
where
\begin{equation}\label{eq:Vdef}
  V=\frac{3(p')^2}{4p^4}-\frac{p''}{2p^3}=\frac{5(q')^2}{16p^6}+\frac{q''}{4p^4}.
\end{equation}
Differentiating with respect to $z$ gives 
\begin{equation}
    \partial_zV=-\frac{15(q')^2}{16p^8}-\frac{q''}{2p^6}
\end{equation}
The derivative estimates above now give 
\begin{equation}
    V\in \Sh_{2}^{r-2},\qquad\partial_z V\in \Th_{2}^{r-2}.
\end{equation}

For a function $\omega$ bounded away from zero, define
\begin{equation}
  \mathcal A(\omega)
  =\frac34\left(\frac{D\omega}{\omega}\right)^2
  -\frac12\frac{D^2\omega}{\omega}.
\end{equation}
Starting with $\omega_0=1$, define recursively
\begin{equation}
  \omega_{n+1}
  =\left(1+V+\mathcal A(\omega_n)\right)^{1/2},
  \qquad n\ge0,
\end{equation}
where we choose the branch tending to $1$ as $k\to\infty$. Since $\sup_{x\geq0}h(x,z)\rightarrow0$ as $\Re z\rightarrow\infty$, the recursion is well defined for $z$ with sufficiently large real part.

\begin{proposition}
For $1\leq n\leq \frac{r-2}{2}$,
\begin{equation}\label{eq:omega-difference}
  \omega_{n+1}-\omega_{n}\in\Sh_{2n+2}^{r-2n-2}.
\end{equation}
Consequently,
\begin{equation}\label{eq:omega-bounds}
  \omega_n\in1+\Sh_2^{r-2n},
  \qquad
  \abs{\omega_n}+\abs{\omega_n^{-1}}\le C_n.
\end{equation}
If we define the defect by
\begin{equation}
  F_n=1+V-\omega_n^2+\mathcal A(\omega_n),
\end{equation}
then
\begin{equation}
  F_n\in\Sh_{2n+2}^{r-2n-2}.
\end{equation}
Furthermore,
\begin{equation}\label{eq:Pz}
  \partial_z\omega_n\in \Th_2^{r-2n},\qquad \partial_z(p\omega_n) =\frac1{2p}\bigl(1+\Oh(h^2)\bigr).
\end{equation}
\end{proposition}

\begin{proof}
We argue by induction. For $n=0$, the definitions give
\begin{equation}
    \omega_0=1\in 1+\Sh_2^{r-2},\qquad F_0=V\in \Sh_2^{r-2},\qquad\partial_z(p\omega_0)=\frac{1}{2p}.
\end{equation}
To treat $\omega_1$, differentiate $\omega_1^2=1+V\in 1+\Sh_2^{r-2}$ and use induction on the derivative order $j$. Assuming the bound $|D^l\omega_1|\leq C_lh^{2+l}$ for the lower positive orders $l<j$, we obtain
\begin{equation}
|2w_1D^jw_1|=\abs{D^jV-\sum_{l=1}^{j-1}\binom jl(D^lw_1)(D^{j-l}w_1)}\leq C(h^{j+2}+h^{j+4}),
\end{equation}
Since $|2\omega_1|\geq1$, this estimate gives the required bound at order $j$. Together with the zeroth-order estimate, it yields $\omega_1-1\in\Sh_2^{r-2}$.

Suppose that $\omega_n-\omega_{n-1}\in\Sh_{2n}^{r-2n}$ and that both functions belong to $1+\Sh_2^{r-2n}$. In the difference of their $\mathcal A$-expressions, each term gains two powers of $h$ from the scaled derivatives. By \eqref{eq:relationship between S,T},
\begin{equation}
  \mathcal A(\omega_{n})-\mathcal A(\omega_{n-1})\in\Sh_{2n+2}^{r-2n-2}.
\end{equation}
Applying the recursion at index $n$ to the difference in \eqref{eq:omega-difference} gives
\begin{equation}
  \omega_{n+1}^2-\omega_n^2=\mathcal A(\omega_n)-\mathcal A(\omega_{n-1})\in\Sh_{2n+2}^{r-2n-2}.
\end{equation}
Since $\omega_{n+1}+\omega_n$ is bounded away from zero, division gives
\begin{equation}
  \omega_{n+1}-\omega_n\in\Sh_{2n+2}^{r-2n-2}.
\end{equation}
This completes the induction.

The recursion also gives
\begin{equation}
  F_n=\mathcal A(\omega_n)-\mathcal A(\omega_{n-1}),
\end{equation}
so $F_n\in\Sh_{2n+2}^{r-2n-2}$.

Finally, differentiating the recursion for $w_n^2$ with respect to $z$ gives
\begin{equation}
    \partial_z\omega_n=\frac{\partial_z V+\partial_z \mathcal{A}(\omega_{n-1})}{2\omega_n}
\end{equation}
The commutator identity $\partial_zD=D\partial_z-\frac{1}{2p^2}D$ gives
\begin{equation}\left\{\begin{aligned}
&\partial_z\left(\frac{D\omega_n}{\omega_n}\right)=\frac{\partial_z(D\omega_n)}{\omega_n}-\frac{D\omega_n}{\omega_n^2}\partial_z\omega_n\in\mathcal \Th_3^{r-2n-1},\\
&\partial_z\left(\frac{D^2\omega_n}{\omega_n}\right)=\frac{\partial_z(D^2\omega_n)}{\omega_n}-\frac{D^2\omega_n}{\omega_n^2}\partial_z\omega_n\in\mathcal \Th_4^{r-2n-2}.\end{aligned}\right.
\end{equation}
Together with \eqref{eq:omega-bounds}, these bounds give
\begin{equation}
  \partial_z\omega_n\in\Th_{2}^{r-2n}.
\end{equation}
Therefore,
\begin{equation}
  \partial_z(p\omega_n)
  =\frac{\omega_n}{2p}+p\partial_z\omega_n
  =\frac1{2p}\bigl(1+\Oh(h^2)\bigr).
\end{equation}
\end{proof}

Define an approximate solution of the transformed Riccati equation by
\begin{equation}\label{eq:Rndef}
  R_n=\ip\omega_n-\frac{D\omega_n}{2\omega_n}.
\end{equation}
A direct calculation gives
\begin{equation}
  D\left(\ip\omega-\frac{D\omega}{2\omega}\right)
  +\left(\ip\omega-\frac{D\omega}{2\omega}\right)^2
  =-\omega^2+\mathcal A(\omega),
\end{equation}
and hence
\begin{equation}
  DR_n+R_n^2+1+V=F_n.
\end{equation}

\subsubsection{A Volterra estimate for the Riccati remainder}

Let $R$ solve the exact equation
\begin{equation}\label{eq:Rexact}
  DR+R^2+1+V=0,
\end{equation}
and write
\begin{equation}\label{eq:Rdelta}
  R=R_n+\delta.
\end{equation}
Then
\begin{equation}\label{eq:delta-D}
  D\delta+2R_n\delta+\delta^2=-F_n,
\end{equation}
or equivalently,
\begin{equation}\label{eq:delta-x}
  \delta'+2pR_n\delta+p\delta^2=-pF_n.
\end{equation}
Introduce the variable $w$ by setting 
\begin{equation}
    \delta=-\frac{2i\omega_nw}{1+w}
\end{equation}
For $|w|\leq\frac{1}{2}$, this defines $\delta$ and gives $|\delta|\leq C|w|$.

In terms of $w$, equation \eqref{eq:delta-x} becomes
\begin{equation}\label{eq:equation of w}
    w'+2ip\omega_nw=b_n(1+w)^2
\end{equation}
where $b_n=\frac{pF_n}{2i\omega_n}$. Taking $2n+2=r$ and using $F_n\in \Sh_{2n+2}^{r-2n-2}$, we obtain
\begin{equation}
    |b_n|\leq C|p|h^r=C|p|^{1-r}\rho^r\leq C|p|^{1-r+2ra}
\end{equation}
Since $a$ satisfies $r(1-2a)\geq M+N+2>2$ and $|p|\geq C|z|^{1/2}$, it follows that
\begin{equation}\label{eq:estimate of b_n}
    |b_n|\leq C|p|^{-1}|z|^{\frac{1}{2}(2-r(1-2a))}\rightarrow0
\end{equation}
Set $S(x,z)=\int_{0}^{x}{|p(t,z)|^{-1}}\dd t$ and define 
\begin{equation}
    Tw(x,z)=-\int_x^\infty\exp\left(2i\int_x^t p(s,z)\omega_n(s,z)\dd s\right)b_n(t,z)(1+w(t,z))^2\dd t.
\end{equation}
We seek a solution of equation~\ref{eq:equation of w} as a fixed point $Tw=w$. To show that $T$ is a contraction, we first estimate its integral kernel.

If $\Im z=0$, both $p$ and $\omega_n$ are real, so $\Im(p(x,\Re z)\omega_n(x,\Re z))=0$.

Now suppose that $0<\Im z<\Re z$. By \eqref{eq:Pz}, increasing $k_0$ if necessary, we obtain
\begin{equation}
  \operatorname{Re}\partial_z(p\omega_n)(x,z)
  \ge \frac{c}{2|p(x,z)|},
  \qquad 0\le\tau\le1.
\end{equation}
Since $p\omega_n$ is real at $z=k^2$, integration in the imaginary direction gives
\begin{equation}
    \begin{aligned}
  \Im (p\omega_n)(x,z)
  &=\int_0^{\Im z}
  \Re\partial_z(p\omega_n)(x,\Re z+\ip\tau)\dd\tau\geq \frac{c\Im z}{2|p(x,z)|}.
\end{aligned}\end{equation}
Combining this with \eqref{eq:estimate of b_n} gives
\begin{equation}\begin{aligned}\label{eq:estimate for kernel}
    \int_x^\infty\abs{\exp\left(2i\int_x^t p(s,z)\omega_n(s,z)\dd s\right)b_n(t,z)}\dd t&\leq C|z|^{\frac{1}{2}(2-r(1-2a))}\int_{x}^{\infty}|p(t,z)|^{-1}e^{-c\Im z(S(t,z)-S(x,z)}\dd t\\
    &\leq C|z|^{\frac{1}{2}(2-r(1-2a))}\int_{0}^{\infty}e^{-c\Im zS}\dd S\leq C' |z|^{\frac{1}{2}(2-r(1-2a))}(1+\frac{1}{\Im z})
\end{aligned}\end{equation}
after the change of variable $S=S(x,z)$. For $\eta>c_0|z|^{-\frac{M}{2}}$, it follows that
\begin{equation}
    \frac{|z|^{\frac{1}{2}(2-r(1-2a))}}{\Im z}\leq c_0^{-1}|z|^{\frac{-N}{2}}.
\end{equation}
Consider $T$ on $C_b[0,\infty)$ and let $B=\{w:\norm{w}_{\infty}\leq\frac{1}{2}\}$. Then 
\begin{equation}
    \norm{(1+w)^2}
    _{\infty}\leq\frac{9}{4},\qquad \norm{(1+w_1)^2-(1+w_2)^2}_{\infty}\leq 3\norm{w_1-w_2}_{\infty}
\end{equation}
Choose $k_0$ so large that the upper bound in \eqref{eq:estimate for kernel} is less than $\frac{2}{9}$. Then
\begin{equation}
    |Tw(x)|\leq \frac{1}{2},\qquad|Tw_1(x)-Tw_2(x)|\leq\frac{2}{3}|w_1(x)-w_2(x)|
\end{equation}
Hence $T:B\rightarrow B$ is a contraction. Its unique fixed point $w$ solves \eqref{eq:equation of w} and satisfies
\begin{equation}\label{eq:estimate of delta}
    \norm{w}_\infty\leq C|z|^{\frac{-N}{2}},\qquad \norm{\delta}_\infty\leq C'\norm{w}_\infty\leq C'' |z|^{\frac{-N}{2}}
\end{equation}

\subsubsection{Identification of the Weyl solution}
Fix $z\in \mathcal{D}\cap \{0<\Im z<\Re z\}$ and choose $n$ so that $r=2n+2$. Set $R=R_n+\delta$ and define $u$ by
\begin{equation}\label{eq:u-definition}
  \frac{Du}{u}=R,
\end{equation}
and set $y=p^{-1/2}u$. Then $y$ solves the original equation \eqref{eq:schrodinger}. By \eqref{eq:Rndef},
\begin{equation}\label{eq:y-formula}
  y(x)=Cp(x)^{-1/2}\omega_n(x)^{-1/2}\exp\left(i\int_0^xp(s)\omega_n(s)\dd s+\int_0^xp(s)\delta(s)\dd s\right).
\end{equation}
Using the relation between $\delta$ and $w$, we obtain 
\begin{equation}
    y^2(x)=\frac{C(1+w(x))^2}{\omega_n(x)}\frac{1}{p(x)}\exp\left(2\int_0^x (ip(s)\omega_n(s)-b_n(s)(1+w(s)))\dd s\right).
\end{equation}
Estimate \eqref{eq:estimate of b_n} gives
\begin{equation}
    \left|\int_0^xb_n(1+w)\dd t\right|\leq C|z|^{\frac{1}{2}(2-r(1-2a))}S(x,z).
\end{equation}
The imaginary part of the phase satisfies 
\begin{equation}
    \Im\int_0^xp\omega_n\dd t\geq c\Im z S(x,z).
\end{equation}
Since $\Im z\geq C |z|^{-\frac{M}{2}}\geq C|z|^{\frac{1}{2}(N+2-r(1-2a))}$, we may choose $k_0$ sufficiently large that \begin{equation}
    \Re\left(2\int_0^x (ip(s)\omega_n(s)-b_n(s)(1+w(s)))\dd s\right)\leq -c\Im z S(x,z)
\end{equation}
which yields
\begin{equation}
    \int_{0}^{\infty}|y^2(x)|\dd x\leq\int_{0}^{\infty}\frac{C''}{|p(x,z)|}e^{-c\Im z S(x,z)}\dd x\leq\frac{C''}{c\Im z}<\infty
\end{equation}
for every $z\in \mathcal{D}\cap \{0<\Im z<\Re z\}$. By the limit-point hypothesis, the square-integrable solution is unique up to a constant factor. Thus $y$ is the Weyl solution of \eqref{eq:schrodinger}, and
\begin{equation}\label{eq:m-logderivative}
  m(z)=\frac{y'(0,z)}{y(0,z)}.
\end{equation}

Set
\begin{equation}\label{eq:r-def}
  m_+(x,z)=\frac{y'(x,z)}{y(x,z)}.
\end{equation}
Since $y=p^{-1/2}u$ and $Du/u=R$,
\begin{equation}\label{eq:r-from-R}
  m_+=-\frac{p'}{2p}+pR.
\end{equation}
Define
\begin{equation}\label{eq:r-app}
  r_n^{\mathrm{app}}
  =-\frac{p'}{2p}
  +p\left(\ip\omega_n-\frac{D\omega_n}{2\omega_n}\right).
\end{equation}
Then
\begin{equation}\label{eq:r-error-exact}
  m_+-r_n^{\mathrm{app}}=p\delta.
\end{equation}
At $x=0$, estimate \eqref{eq:estimate of delta} gives 
\begin{equation}\label{eq:r0-small-error}
    \abs{m_+(0,z)-r_n^{\mathrm{app}}(0,z)}\leq C'''|z|^{\frac{1}{2}(-N+1)}
\end{equation}

At $x=0$, the derivatives of $q$ that enter the construction are finite, and
\begin{equation}\label{eq:p0-expansion}
  p(0,z)
  =\sqrt{z}\left(1-\frac{q(0)}{z}\right)^{1/2}
\end{equation}
has a complete asymptotic expansion in powers of $\sqrt{z}^{-1}$. The expression $r_n^{\mathrm{app}}(0,z)$ involves only $p$, finitely many of its $x$-derivatives, algebraic operations, and square roots. It therefore also has an expansion in powers of $\sqrt{z}^{-1}$.

The approximate logarithmic derivative satisfies the exact defect identity
\begin{equation}\label{eq:original-defect}
  (r_n^{\mathrm{app}})'
  +(r_n^{\mathrm{app}})^2-q+z
  =p^2F_n.
\end{equation}
On each fixed compact interval, $h=\Oh(|z|^{-\frac{1}{2}})$ and $p=\Oh(|z|^{\frac{1}{2}})$. Since $r\geq N+2$, the estimate
\begin{equation}\label{eq:defect-compact}
  p^2F_n=\Oh\left(|z|^{\frac{-N}{2}}\right)
\end{equation}
holds uniformly, together with the finitely many $x$-derivatives needed below. Thus, on compact intervals,
\begin{equation}\label{eq:local-expansion}
  r_n^{\mathrm{app}}(x,z)
  =\sqrt{-z}+\sum_{j=1}^{N-2}\frac{c_j(x)}{(2\sqrt{-z})^j}
  +\Oh\left(|z|^{\frac{1}{2}(-N+1)}\right).
\end{equation}

Evaluating at $x=0$ and using \eqref{eq:m-logderivative} and \eqref{eq:r0-small-error} gives, uniformly for $z\in \mathcal{D}\cap \{0<\Im z<\Re z\}$,
\begin{equation}
  m_+(z)
  =\sqrt{-z}+\sum_{j=1}^{N-2}\frac{c_j(0)}{(2\sqrt{-z})^j}
  +\Oh\left(|z|^{\frac{1}{2}(-N+1)}\right).
\end{equation}
The relations among the coefficients and the expansion in the sector $\mathcal{D}\cap\{\Im z>a\Re z\geq C_0+1\}$ are given in \cite{DanielyanLevitan1991}.

Applying the same argument to the left half-line Weyl function $m_-$ and to the other branch of $\mathcal{D}$ gives the Weyl expansions required in Theorem~\ref{Thm: theorem3}.

\subsection{Exponential decay of the spectral measure}
We use the following estimate of Gesztesy and Simon (\cite[Theorem~6.6]{GesztesySimon2000}). 
\begin{lemma}
    Suppose that $q \in L^1(0,a)$ for every $a > 0$ and that $+\infty$ is a limit-point endpoint. Then, for every $\delta\in(0,1)$ and $\alpha_0 > 0$,
    \begin{equation}
        \alpha_0\delta\int_{-\infty}^{0}e^{2(1-\delta)\alpha_0\sqrt{-\lambda}}\sigma_+(d\lambda)\leq C_1(1+\alpha_0)+C_2(1+E(\alpha_0)^2)e^{2(1+\alpha_0)\sqrt{E(\alpha_0)}}
    \end{equation}
    where
    \begin{equation}
        E(\alpha_0)=-\inf_{\varphi}\left\{\int_{0}^{1+\alpha_0}|\varphi'(x)|^2+q(x)|\varphi(x)|^2 dx:\varphi\in C_{0}^{\infty}((0,\alpha_0+1)),\|\varphi\|_{L^2((0,\alpha_0+1))}\leq1\right\}.
    \end{equation}
    The constants $C_1,C_2$ depend only on $L=\int_0^1|q(x)|dx$.
\end{lemma}
This estimate gives the following exponential moment bound.
\begin{lemma}
    Suppose that $q(x)\geq-Cx^\alpha$ for $x>x_0$, with $x_0$ sufficiently large. Then, for sufficiently small $\varepsilon>0$,
    \begin{equation}
        \int_{-\infty}^{0}e^{\varepsilon|\lambda|^{\frac{1}{2}+\frac{1}{\alpha}}}\sigma_+(\dd\lambda)<\infty
    \end{equation}
\end{lemma}
\begin{proof}
     Set $K_0=\int_{0}^{x_0}q_-(x)\dd x<\infty$, where $q_-(x)=\max\{-q(x),0\}$. For $\varphi\in C^{\infty}_0$, we have 
    \begin{equation}
        \norm{\varphi}_{\infty}^2\leq2\left(\int_{0}^{1+\alpha_0} \abs{\varphi'(s)\overline{\varphi(s)}}\dd s\right)\leq 2\norm{\varphi'}_{2}\norm{\varphi}_2
    \end{equation}
    where $\norm{\cdot}_2=\norm{\cdot}_{L^2((0,\alpha_0+1))}$. Splitting the integral over $(0,x_0)$ and $(x_0,1+\alpha_0)$ gives
    \begin{equation}
        \int_{0}^{1+\alpha_0}-q(x)|\varphi(x)|^2 \dd x\leq K_0\norm{\varphi}_\infty^2+C(1+\alpha_0)^\alpha\norm{\varphi}_2^2\leq 2K_0\norm{\varphi'}_{2}\norm{\varphi}_2+C(1+\alpha_0)^\alpha\norm{\varphi}_2^2.
    \end{equation}
    Completing the square gives
\begin{equation}\begin{aligned}
    -\int_{0}^{1+\alpha_0}|\varphi'(x)|^2+q(x)|\varphi(x)|^2 \dd x&\leq-\norm{\varphi'}_{2}^{2}+ 2K_0\norm{\varphi'}_{2}\norm{\varphi}_2+C(1+\alpha_0)^\alpha\norm{\varphi}_2^2\\
    &\leq -(\norm{\varphi'}_2
    -K_0\norm{\varphi}_2)^2+(K_0^2+C(1+\alpha_0)^\alpha)\norm{\varphi}_2^2
    \\
    &\leq  K_0^2+C(1+\alpha_0)^\alpha
\end{aligned}
\end{equation}
for every $\varphi$ with $\|\varphi\|_{2}\leq1$. Taking the supremum gives
\begin{equation}
    E(\alpha_0)=\sup_{\varphi}\left\{-\int_{0}^{1+\alpha_0}|\varphi'(x)|^2+q(x)|\varphi(x)|^2 \dd x\right\}\leq K_0^2+C(1+\alpha_0)^\alpha
\end{equation}

Taking $\delta=\frac{1}{2}$ in the estimate of Gesztesy and Simon gives 
\begin{equation}
    \int_{-\infty}^{0}e^{B|\lambda|^{\frac{1}{2}}}\sigma_+(d\lambda)\leq D_1e^{D_2 B^{1+\frac{\alpha}{2}}}
\end{equation}
where $B=2(1-\delta)\alpha_0=\alpha_0$. For each fixed $\varepsilon_1>0$, we may take $D_2=2\sqrt{C+\varepsilon_1}$, with a constant $D_1$ depending on $C,L,K_0,\alpha,\varepsilon_1$.

Decomposing the negative half-line into the intervals between consecutive negative squares gives
\begin{equation}
    \begin{aligned}
        \int_{-\infty}^{0}e^{\varepsilon|\lambda|^{t}}\sigma_+(\dd\lambda)&=\sum_{n=0}^{\infty}\int_{-(n+1)^2}^{-n^2}e^{\varepsilon|\lambda|^{t}}\sigma_+(\dd\lambda)\leq\sum_{n=0}^{\infty}e^{\varepsilon(n+1)^{2t}}\int_{-(n+1)^2}^{-n^2}\sigma_+(\dd\lambda)\\
        &\leq\sum_{n=0}^{\infty}e^{\varepsilon(n+1)^{2t}}\int_{-(n+1)^2}^{-n^2}e^{B|\lambda|^{\frac{1}{2}}}e^{-Bn}\sigma_+(\dd\lambda)\\
        &\leq\sum_{n=0}^{\infty}e^{\varepsilon(n+1)^{2t}-Bn} D_1e^{D_2B^{1+\frac{1}{2}\alpha}}\\
        &\leq D_1\sum_{n=0}^{\infty}\exp\left({\varepsilon(n+1)^{2t}-Bn+{D_2B^{1+\frac{1}{2}\alpha}}}\right).
    \end{aligned}
\end{equation}
On each interval $(-(n+1)^2,-n^2]$, choose $B=\left(\frac{n}{2D_2}\right)^{\frac{2}{\alpha}}$. Then
\begin{equation}
     \int_{-\infty}^{0}e^{\varepsilon|\lambda|^{t}}\sigma_+(d\lambda)\leq D_1\sum_{n=0}^{\infty}\exp\left(\varepsilon(n+1)^{2t}-\frac{1}{2}\left(\frac{1}{2D_2}\right)^{\frac{2}{\alpha}}n^{\left(1+\frac{2}{\alpha}\right)} \right).
\end{equation}

The series on the right converges, and hence 
\begin{equation}
    \int_{-\infty}^{0}e^{\varepsilon|\lambda|^{\frac{1}{2}+\frac{1}{\alpha}}}\sigma_+(d\lambda)<\infty
\end{equation}
for $t=\frac{1}{2}+\frac{1}{\alpha}$ and every $0<\varepsilon\leq 2^{-\frac{4}{\alpha}-1}C^{-\frac{1}{\alpha}}$.

\end{proof}
Combining this estimate with the Weyl expansion established under the hypotheses of Lemma~\ref{lem: expansion} proves Theorem~\ref{Thm: theorem3}.

\end{document}